\documentclass[openright,11pt,two0side,a4paper,leqno]{article}
\usepackage[latin1]{inputenc} 
\usepackage{amsmath,amssymb} 
\usepackage{theorem}

\newlength{\larg}
\newcommand{\bd}[1]{\boldsymbol #1}
\newcommand{\ds}{\displaystyle}

\newcommand{\pv}{\noindent \emph{Proof. }}
\newcommand{\cqfd}{\phantom{m} \hfill $\Box$}
\newcommand{\fini}{\phantom{m} \hfill $\diamond$}
\newenvironment{proof}{\pv}{\cqfd \\ }

\newcommand{\eg}{\emph{e.g. }}
\newcommand{\resp}{\emph{resp. }}

\newcommand{\barre}{\overline{ ^{\ \: }} \ }

\newcommand{\set}[1]{ \{ #1 \}}

\newcommand{\mA}{{\mathbb A}}
\newcommand{\mC}{{\mathbb C}}
\newcommand{\mN}{{\mathbb N}}
\newcommand{\mQ}{{\mathbb Q}}
\newcommand{\mR}{{\mathbb R}}
\newcommand{\mZ}{{\mathbb Z}}

\newcommand{\interv}[2]{[\! [ #1  ;  #2 ] \!]}

\newcommand{\wt}{\mathrm{wt}}
\newcommand{\wtpt}{\dot{\mathrm{wt}}}

\newcommand{\eiup}{\tilde{e_i}^{\scriptsize \mbox{\rm up}}}
\newcommand{\eilow}{\tilde{e_i}^{\scriptsize \mbox{\rm low}}}
\newcommand{\fiup}{\tilde{f_i}^{\scriptsize \mbox{\rm up}}}
\newcommand{\filow}{\tilde{f_i}^{\scriptsize \mbox{\rm low}}}

\newcommand{\gsl}[1]{{\mathfrak{sl}}_{#1}}
\newcommand{\hsl}[1]{{\widehat{\mathfrak{sl}}}_{#1}}
\newcommand{\Uq}{U_q(\hsl{n})}
\newcommand{\Up}{U_p(\hsl{l})}
\newcommand{\Uqprime}{U_q'(\hsl{n})}
\newcommand{\Upprime}{U_p'(\hsl{l})}

\newcommand{\Fq}[1][\bd{s}_l]{\boldsymbol{\operatorname{F}}_q[#1]}
\newcommand{\Fp}[1][\bd{s}_n]{\boldsymbol{\operatorname{F}}_p[#1]^{\bullet}}

\newcommand{\res}{\operatorname{res}}

\newcommand{\zconnected}[1]{\ _{\stackrel{_{- \hspace{-3mm} - \hspace{-3mm} - \hspace{-3mm} - }}{#1}} \ }
\newcommand{\Zconnected}[1]{\ _{\stackrel{\equiv}{#1}} \ }

\newtheorem{thm}{Theorem}[section]
\newtheorem{lemma}[thm]{Lemma}
\newtheorem{prop}[thm]{Proposition}
\newtheorem{cor}[thm]{Corollary}

{\theorembodyfont{\rmfamily} \newtheorem{remark}[thm]{Remark}}   
{\theorembodyfont{\rmfamily} \newtheorem{example}[thm]{Example}} 
{\theorembodyfont{\rmfamily} \newtheorem{definition}[thm]{Definition}}
{\theorembodyfont{\rmfamily} \newtheorem{notation}[thm]{Notation}}

\title{Canonical bases of higher-level $q$-deformed Fock spaces}
\author{Xavier YVONNE}
\date{}
\begin{document}

\maketitle

\begin{abstract} We show that the transition matrices between the standard and the canonical bases of infinitely many
weight subspaces of the higher-level $q$-deformed Fock spaces are equal. 
\end{abstract}

\section{Introduction}
The $q$-deformed higher-level Fock spaces were introduced in \cite{JMMO} in order to compute the crystal graph of
any irreducible integrable representation of level $l \geq 1$ of $\Uq$. More precisely, the Fock representation $\Fq$ depends
on a parameter $\bd{s}_l=(s_1,\ldots,s_l) \in \mZ^l$ called multi-charge. It contains as a submodule the irreducible 
integrable $\Uq$-module with highest weight $\Lambda_{s_1}+\cdots+\Lambda_{s_l}$. The representation $\Fq$ is a generalization of
the level-one Fock representation of $\Uq$ (\cite{H,MM}, see also \cite{LT1,LT2}). 

The canonical bases are bases of the Fock representations that are invariant under a certain involution $\barre$ of
$\Uq$ and that give at $q=0$ and $q=\infty$ the crystal bases. They were constructed for $l=1$ in \cite{LT1,LT2} and for $l \geq 1$ by
Uglov \cite{U}. In \cite{U}, Uglov provides an algorithm for computing these canonical bases. He also gives an expression
of the transition matrices between the standard and the canonical bases in terms of Kazhdan-Lusztig polynomials for affine Hecke algebras of type $A$.

In this paper we prove three theorems.
\begin{enumerate}
\item The first one (Theorem \ref{thm_comparaison_bc_1}) is a generalization to $l \geq 1$
of a result of \cite{LM}. It compares the transition matrices of the canonical bases of some weight subspaces inside a given Fock space $\Fq$.
The weights involved are conjugated under the action of the Weyl group of $\Uq$. This action leads to bijections $\sigma_i$ that can be described in a combinatorial
way by adding/removing as many $i$-nodes as possible to the $l$-multi-partitions indexing the canonical bases. These bijections are generalizations of the Scopes bijections
introduced in \cite{S} in order to study, when $n=p$ is a prime number, the $p$-blocks of symmetric groups of a given defect.  
\item In a dual manner, our second result (Theorem \ref{thm_comparaison_bc_2}) gives some sufficient conditions on multi-charges
$\bd{s}_l$ and $\bd{t}_l$ with given residues modulo $n$ that ensure that the transition matrices of the canonical bases of some weight subspaces
of $\Fq$ and $\Fq[\bd{t}_l]$ coincide.
\item Our third result (Theorem \ref{thm_comparaison_bc_cas_dominant}) is an application of Theorem \ref{thm_comparaison_bc_2} to the case when
the multi-charges $\bd{s}_l=(s_1,\ldots,s_l)$ and $\bd{t}_l=(t_1,\ldots,t_l)$ are dominant, that is $s_1 \gg \cdots \gg s_l$ and $t_1 \gg \cdots \gg t_l$.
It shows that the transition matrices of the canonical bases of the Fock spaces $\Fq$ stabilize when $\bd{s}_l$ becomes
dominant (with a given sequence of residues modulo $n$). This supports the following conjecture (see \cite{Y2}). We
conjecture that if $\bd{s}_l=(s_1,\ldots,s_l)$ is dominant, then the transition matrix of the homogeneous component of degree $m$ of
the canonical basis of the Fock space $\Fq$ is equal to the decomposition matrix of the cyclotomic $v$-Schur algebra
$\mathcal{S}_{\mC,m}(\zeta;\zeta^{s_1},\ldots,\zeta^{s_l})$ of \cite{DJM}, where $\zeta$ is a complex primitive $n$-th root
of unity. This conjecture generalizes both Ariki's theorem for Ariki-Koike algebras (see \cite{A}) and a result of
Varagnolo and Vasserot (see \cite{VV}) which relates the canonical basis of the level-one Fock space and the decomposition
matrix of $v$-Schur algebras with parameter a complex $n$-th root of unity. 
\end{enumerate}

\noindent \emph{Acknowledgments.} I thank Bernard Leclerc for his advice and his reading of the preliminary 
versions of this paper.\\

\noindent \emph{Notation.} Let $\mN$ (\resp $\mN^{*}$) denote the set of nonnegative (\resp positive) integers, and for $a,\,b \in \mR$ denote by 
$\interv{a}{b}$ the discrete interval $[a  ,  b] \cap \mZ$. For $X \subset \mR$, $t \in \mR$, $N \in \mN^*$, put 
\begin{equation} \label{eq_Z_l_s}
X^N(t):=\set{(s_1,\ldots,s_N) \in X^N \mid s_1+\cdots+s_N=t}.
\end{equation}
Throughout this article, we fix 3 integers $n$, $l \geq 1$ and $s \in \mZ$. Let $\Pi$
denote the set of all integer partitions. If $(W,S)$ is a Coxeter system, denote by $\ell : W \rightarrow \mN$ the
length function on $W$.

\section{Higher-level $q$-deformed Fock spaces}
In this section, we introduce the higher-level Fock spaces and their canonical bases. We follow here \cite{U}, to which
we refer the reader for more details. All definitions and results given here are due to Uglov. 

\subsection{The quantum algebras $\Uq$ and $\Up$}

In this section, we assume that $n \geq 2$ and $l \geq 2$. Let $\hsl{n}$ be the Kac-Moody algebra of type $A_{n-1}^{(1)}$ defined over the field 
$\mQ$ \cite{Kac}. Let $\mathfrak{h}^*$ be the dual of the Cartan subalgebra of $\hsl{n}$. Let $\Lambda_0,\ldots,\Lambda_{n-1} \in \mathfrak{h}^*$ 
be the fundamental weights, $\alpha_0,\ldots,\alpha_{n-1} \in \mathfrak{h}^*$ be the simple roots and $\delta:=\alpha_0+\cdots+\alpha_{n-1}$ be 
the null root. It will be convenient to extend the index set of the fundamental weights by setting $\Lambda_i:=\Lambda_{i \bmod n}$ for all 
$i \in \mZ$. The simple roots are related to the fundamental weights by
\begin{equation} 
\alpha_i=2 \Lambda_i - \Lambda_{i-1} - \Lambda_{i+1} + \delta_{i,0} \delta \qquad (0 \leq i \leq n-1).
\end{equation} 
For $0 \leq i,j \leq n-1$, let $a_{ij}$ be the coefficient of $\Lambda_j$ in $\alpha_i$. The space
$$\mathfrak{h}^* = \bigoplus_{i=0}^{n-1}{\mQ\,\Lambda_i} \; \oplus \; \mQ\,\delta = \bigoplus_{i=0}^{n-1}{\mQ\,\alpha_i} \; \oplus \; \mQ\,\Lambda_0$$ 
is equipped with a non-degenerate bilinear symmetric form $(.\,,\,.)$ defined by
\begin{equation}
(\alpha_i,\alpha_j)=a_{ij}, \quad (\Lambda_0,\alpha_i) = \delta_{i,0}, \quad 
(\Lambda_0,\Lambda_0)=0 \qquad (0 \leq i,j \leq n-1).
\end{equation} 

Let $\Uq$ be the $q$-deformed universal enveloping algebra of $\hsl{n}$. This is an algebra over $\mQ(q)$ with generators $e_i$, $f_i$, $t_i^{\pm 1}$
$(0 \leq i \leq n-1)$ and $\partial$. Let $\Uqprime$ be the subalgebra of $\Uq$ generated by $e_i$, $f_i$, $t_i^{\pm 1}$ $(0 \leq i \leq n-1)$.
The relations in $\Uqprime$ are standard and will be omitted (see \eg \cite{KMS}). The relations among the degree generator $\partial$ and
the generators of $\Uqprime$ can be found in \cite{U}. If $M$ is a $\Uq$-module, denote by
$\mathcal{P}(M)$ the set of weights of $M$ and let $M \langle w \rangle$ denote the subspace of $M$ of weight $w$.
If $x \in M \langle w \rangle \setminus \set{0}$ is a weight vector, denote by
\begin{equation}  
\wt(x):=w
\end{equation}
the weight of $x$. The \emph{Weyl group} of $\hsl{n}$ (or $\Uq$), denoted by $W_n$, is the subgroup of $\mathrm{GL}(\mathfrak{h}^*)$ generated by the 
simple reflections $\sigma_i$ defined by
\begin{equation}  \label{def_Weyl_group}
\sigma_i(\Lambda)=\Lambda - (\Lambda,\alpha_i)\, \alpha_i \qquad (\Lambda \in \mathfrak{h}^*,\, 0 \leq i \leq n-1).
\end{equation}
Note that $W_n$ is isomorphic to $\widetilde{\mathfrak{S}}_n$, the affine symmetric group which is a Coxeter group of type $A_{n-1}^{(1)}$. \\

We also introduce the algebra $\Up$ with
\begin{equation}
p:=-q^{-1}.
\end{equation}
In order to distinguish the elements related to $\Uq$ from those related to $\Up$, we put dots over the latter. For example, $\dot{e}_i$, 
$\dot{f}_i$, $\dot{t}_i^{\pm 1}$ $(0 \leq i \leq l-1)$ and $\dot{\partial}$ are the generators of $\Up$, $\dot{\alpha}_i$ $(0 \leq i \leq l-1)$ 
are the simple roots for $\Up$, $\dot{W}_l = \langle \dot{\sigma}_0,\ldots,\dot{\sigma}_{l-1} \rangle$ 
is the Weyl group of $\Up$ and so on.

\subsection{The space $\Lambda^s$}
\subsubsection{The vector space $\Lambda^s$ and its standard basis}
Following \cite{U}, we now recall the definition of $\Lambda^s$, the space of (semi-infinite) $q$-wedge products of charge $s$ 
(this space is denoted by $\Lambda^{s+\frac{\infty}{2}}$ in \cite{U}). Let $P(s)$ be the set of the sequences of integers $(k_1,k_2,\ldots)$ such 
that $k_i = s-i+1$ for $i$ large enough. As a vector space over $\mQ(q)$, $\Lambda^s$ is spanned by the \emph{$q$-wedge products}
\begin{equation}
u_{\bd{k}} = u_{k_1} \wedge u_{k_2} \wedge \cdots, \qquad \bd{k}=(k_1,k_2,\ldots) \in P(s),
\end{equation}
with relations given in \cite[Prop. 3.16]{U}. These relations are called \emph{straightening rules}. Using them, one can express a $q$-wedge product
as a linear combination of so-called \emph{ordered $q$-wedge products}, namely $q$-wedge products $u_{\bd{k}}$ with $\bd{k} \in P^{++}(s)$, where
\begin{equation}
P^{++}(s):=\set{ (k_1,k_2,\ldots) \in P(s) \mid k_1>k_2> \cdots}.
\end{equation}
In fact, the ordered $q$-wedge products $\set{u_{\bd{k}} \mid \bd{k} \in P^{++}(s)}$ form a basis of $\Lambda^s$, called the \emph{standard basis}. 
For the sequel, il will be convenient to use different indexations of this basis which we give now.

\begin{itemize}
\item[*] \emph{\underline{Indexation $u_{\bd{k}}$}}. This is the indexation we have just described. \\

\item[*] \emph{\underline{Indexation $\lambda$}}. To the ordered $q$-wedge product $u_{\bd{k}}$ corresponds a partition 
$\lambda=(\lambda_1,\lambda_2,\ldots)$ defined by
\begin{equation}
\lambda_i := k_i-(s+1-i) \qquad (i \geq 1).
\end{equation}
If $u_{\bd{k}}$ and $\lambda$ are related this way, write
\begin{equation}
|\lambda,s \rangle :=u_{\bd{k}}.
\end{equation}

\item[*] \emph{\underline{Indexation $\bd{\lambda}_n$}}. Recall the definition of $\mZ^n(s)$ from (\ref{eq_Z_l_s}).
Uglov constructed a bijection 
\begin{equation}
\tau'_n : \Pi \rightarrow \Pi^n \times \mZ^n(s), \quad \lambda \mapsto (\bd{\lambda}_n,\bd{s}_n)
\end{equation}
(see \cite[Sect. 4.1]{U}, where this map is denoted by $\tau_n^s$). With the notation above, $\bd{\lambda}_n$ is the $n$-quotient of $\lambda$ and 
$\bd{s}_n$ is a variation of the $n$-core of $\lambda$ (see \eg \cite[Ex.8, p.12]{Mac}). Write
\begin{equation}
|\bd{\lambda}_n,\bd{s}_n \rangle^{\bullet} := |\lambda,s \rangle
\end{equation}
if $(\bd{\lambda}_n,\bd{s}_n)=\tau'_n(\lambda)$. Note that this indexation coincides with the indexation $\lambda$ if $n=1$. \\

\item[*] \emph{\underline{Indexation $\bd{\lambda}_l$}}. Uglov constructed a bijection 
\begin{equation}
\tau_l : \Pi \rightarrow \Pi^l \times \mZ^l(s), \quad \lambda \mapsto (\bd{\lambda}_l,\bd{s}_l)
\end{equation}
(see again \cite[Sect. 4.1]{U}, where this map is denoted by $\tau_l^s$). The map $\tau_l$ is a variation of the map $\tau'_n$ defined above. Write
\begin{equation}
|\bd{\lambda}_l,\bd{s}_l \rangle := |\lambda,s \rangle
\end{equation}
if $(\bd{\lambda}_l,\bd{s}_l)=\tau_l(\lambda)$. Note that this indexation coincides with the indexation $\lambda$ if $l=1$. \\
\end{itemize}

\begin{example} 
Take $n=2$ and $l=3$. Then we have
$$\big|(4,3,3,2,1),-1 \big\rangle = \big| \bigl( (3,3),\emptyset \bigr),(-1,0) \big\rangle^{\bullet} = \big| \bigl( (1,1),(1,1),(1) \bigr),(0,0,-1) \big\rangle.$$
\fini
\end{example}  

\subsubsection{3 actions on $\Lambda^s$} \label{paragraphe_3_actions_on_Lambdas}

Following \cite{JMMO},\cite{FLOTW},\cite{U}, the vector space $\Lambda^s$ can be made into an integrable representation of level $l$ of the quantum 
algebra $\Uq$. This representation can be described in a nice way if we use the indexation $\bd{\lambda}_l$. In order to recall the explicit 
formulas, let us first introduce some notation. Fix $\bd{\lambda}_l=(\lambda^{(1)},\ldots,\lambda^{(l)}) \in \Pi^l$ and  
$\bd{s}_l=(s_1,\ldots,s_l) \in \mZ^l$. Identify the multi-partition $\bd{\lambda}_l$ with its Young diagram 
$\set{(i,j,b)\in \mN^*\times\mN^*\times \interv{1}{l} \mid 1 \leq j \leq \lambda^{(b)}_i}$, whose elements are called \emph{nodes} of 
$\bd{\lambda}_l$. For each node $\gamma=(i,j,b)$ of $\bd{\lambda}_l$, define its \emph{residue modulo $n$} by
\begin{equation}
\res_n(\gamma) \ = \ \res_n(\gamma,\bd{s}_l) \ := \ (s_b+j-i) \, \bmod n \; \in \, \mZ/n\mZ \cong \interv{0}{n-1}.
\end{equation}   
If $\res_n(\gamma)=c$, we say that $\gamma$ is a \emph{$c$-node}. If $\bd{\mu}_l \in \Pi^l$ is such that 
$\bd{\mu}_l \supset \bd{\lambda}_l$ and $\gamma:=\bd{\mu}_l \setminus \bd{\lambda}_l$ is a $c$-node of $\bd{\mu}_l$, we say that 
$\gamma$ is a \emph{removable $c$-node of $\bd{\mu}_l$} or that $\gamma$ is an \emph{addable $c$-node of $\bd{\lambda}_l$}. For $0 \leq c \leq n-1$, 
denote by 
\begin{equation}
\mbox{$N_c(\bd{\lambda}_l;\bd{s}_l;n) \quad$ (\resp $A_c(\bd{\lambda}_l;\bd{s}_l;n), \quad$ \resp $R_c(\bd{\lambda}_l;\bd{s}_l;n)$)}
\end{equation}
the number of $c$-nodes (\resp of addable $c$-nodes, \resp of removable $c$-nodes) of $\bd{\lambda}_l$. Put
\begin{equation} 
M_c(\bd{\lambda}_l;\bd{s}_l;n):=A_c(\bd{\lambda}_l;\bd{s}_l;n)-R_c(\bd{\lambda}_l;\bd{s}_l;n).
\end{equation} 
For $\bd{\lambda}_l$, $\bd{\mu}_l \in \Pi^l$, $\bd{s}_l \in \mZ^l$, $c \in \mZ$ and $k \in \mN^{*}$, write
\begin{equation}
\bd{\lambda}_l \stackrel{c:k}{\longrightarrow} \bd{\mu}_l
\end{equation}
if there exists a sequence of $l$-multi-partitions 
$\bd{\nu}_l^{(0)} \subset \bd{\nu}_l^{(1)} \subset \cdots \subset \bd{\nu}_l^{(k)}$ such that
$\bd{\lambda}_l = \bd{\nu}_l^{(0)}$, $\bd{\mu}_l=\bd{\nu}_l^{(k)}$ and for all $1 \leq j \leq k$, $\bd{\nu}_l^{(j)} \setminus \bd{\nu}_l^{(j-1)}$ is an 
addable $c$-node of $\bd{\nu}_l^{(j-1)}$. Given a multi-charge $(s_1,\ldots,s_l)$ and two nodes $\gamma=(i,j,b)$ and $\gamma'=(i',j',b')$, write 
\begin{equation}
\gamma<\gamma'
\end{equation}
if either $s_b+j-i<s_{b'}+j'-i'$ or $s_b+j-i=s_{b'}+j'-i'$ and $b<b'$. This defines a total ordering on the set of the addable and
removable $c$-nodes of a given multi-partition. If $\bd{\lambda}_l \stackrel{c:k}{\longrightarrow} \bd{\mu}_l$, put

\begin{equation}
\begin{array}{rcl}
M_c^{>}(\bd{\lambda}_l;\bd{\mu}_l;\bd{s}_l;n) &=& \ds \sum_{\gamma \in \bd{\mu}_l \setminus \bd{\lambda}_l}
\Bigl( \sharp \set{\beta \in \mN^3 \mid \beta \mbox{ is an addable $c$-node of $\bd{\mu}_l$ and } \beta > \gamma} \\[3mm]
&&\qquad - \sharp \set{\beta \in \mN^3 \mid \beta \mbox{ is a removable $c$-node of $\bd{\lambda}_l$ and } \beta > \gamma} \Bigr), 
\end{array}
\end{equation}
and define similarly $M_c^{<}(\bd{\lambda}_l;\bd{\mu}_l;\bd{s}_l;n)$.

\begin{example} 
Take $\bd{s}_l=(5,0,2,1)$, $\bd{\lambda}_l=\bigl( (5,3,3,1), (3,2), (4,3,1), (2,2,2,1) \bigr)$, $n=3$ and $c=0$. Then we have 
$$N_c(\bd{\lambda}_l;\bd{s}_l;n)=11, \quad A_c(\bd{\lambda}_l;\bd{s}_l;n) = R_c(\bd{\lambda}_l;\bd{s}_l;n) =5 \quad \mbox{and} 
\quad M_c(\bd{\lambda}_l;\bd{s}_l;n)=0.$$
The addable $c$-nodes of $\bd{\lambda}_l$ are $(5,1,4)$, $(4,2,1)$, $(1,4,2)$, $(1,3,4)$ and $(1,5,3)$. 
The removable $c$-nodes of $\bd{\lambda}_l$ are $(2,2,2)$, $(3,1,3)$, $(3,2,4)$, $(2,3,3)$ and $(1,5,1)$. 
The list of all these nodes arranged with respect to the ordering described above is
$$\footnotesize \begin{array}{r} (5,1,4) < (2,2,2) < (3,1,3) < (3,2,4) < (4,2,1) < (1,4,2) < (2,3,3) < (1,3,4) < (1,5,3) < (1,5,1). \end{array}$$
Take also $\bd{\mu}_l=\bigl( (5,3,3,1), (3,2), (5,3,1), (2,2,2,1) \bigr)$, so that $\bd{\mu}_l \setminus \bd{\lambda}_l=\set{(1,5,3)}$
is a single $c$-node. Then $M_c^{>}(\bd{\lambda}_l;\bd{\mu}_l;\bd{s}_l;n)=0-1=-1$ and $M_c^{<}(\bd{\lambda}_l;\bd{\mu}_l;\bd{s}_l;n)=4-4=0$. \fini 
\end{example}

For $a \in \mZ$, denote by $a \bmod n \in \set{0,\ldots,n-1}$ the residue of $a$ modulo $n$. Finally, for
$\bd{s}_l = (s_1,\ldots,s_l) \in \mZ^l$, define
\begin{equation}
\Delta(\bd{s}_l,n):= \ds \frac{1}{2} \sum_{b=1}^{l}{\Bigl( \frac{s_b^2}{n}-s_b \Bigr) - 
\Bigl(\frac{(s_b \bmod n)^2}{n}-(s_b \bmod n)\Bigr)}.
\end{equation}

Now we can state the following result.

\begin{thm}[\cite{JMMO},\cite{FLOTW},\cite{U}] \label{thm_action_Uq_sur_Lambda_s} The following formulas define on $\Lambda^s$ a structure of
an integrable representation of level $l$ of the quantum algebra $\Uq$. 
\begin{equation}
e_i.|\bd{\nu}_l, \bd{s}_l \rangle = \ds \sum_{\bd{\lambda}_l \stackrel{i:1}{\longrightarrow} \bd{\nu}_l}
{q^{-M_i^{<}(\bd{\lambda}_l;\bd{\nu}_l;\bd{s}_l;n)}\,|\bd{\lambda}_l,\bd{s}_l \rangle},
\end{equation}
\begin{equation}
f_i.|\bd{\nu}_l,\bd{s}_l \rangle =
\ds \sum_{\bd{\nu}_l \stackrel{i:1}{\longrightarrow} \bd{\mu}_l} 
{q^{M_i^{>}(\bd{\nu}_l;\bd{\mu}_l;\bd{s}_l;n)}\, |\bd{\mu}_l,\bd{s}_l \rangle},
\end{equation}
\begin{equation}
t_i.|\bd{\nu}_l,\bd{s}_l \rangle = q^{M_i(\bd{\nu}_l;\bd{s}_l;n)}\, |\bd{\nu}_l,\bd{s}_l \rangle,
\end{equation}
\begin{equation}
\partial.|\bd{\nu}_l,\bd{s}_l \rangle = -\bigl( \Delta(\bd{s}_l,n) + N_0(\bd{\nu}_l;\bd{s}_l;n) \bigr) \,
|\bd{\nu}_l,\bd{s}_l \rangle. 
\end{equation} \cqfd
\end{thm}
Note that these formulas do no involve the straightening of $q$-wedge products. They are therefore handy to use for computations. \\ 

In a completely similar way, $\Lambda^s$ can be made into an integrable representation of level $n$ of the quantum algebra $\Up$.
This action can be described using the indexation $\bd{\lambda}_n$. Namely, we have (with obvious notation) the following result.

\begin{thm}[\cite{JMMO},\cite{FLOTW},\cite{U}] \label{thm_action_Up_sur_Lambda_s} The following formulas define on $\Lambda^s$ a structure of
an integrable representation of level $n$ of the quantum algebra $\Up$. 
\begin{equation}
\dot{e}_i.|\bd{\nu}_n, \bd{s}_n \rangle^{\bullet} =
\ds \sum_{\bd{\lambda}_n \stackrel{i:1}{\longrightarrow} \bd{\nu}_n}
{p^{-M_i^{<}(\bd{\lambda}_n;\bd{\nu}_n;\bd{s}_n;l)} \,|\bd{\lambda}_n,\bd{s}_n \rangle^{\bullet}},
\end{equation}
\begin{equation}
\dot{f}_i.|\bd{\nu}_n,\bd{s}_n \rangle^{\bullet} =
\ds \sum_{\bd{\nu}_n \stackrel{i:1}{\longrightarrow} \bd{\mu}_n} 
{p^{M_i^{>}(\bd{\nu}_n;\bd{\mu}_n;\bd{s}_n;l)}\, |\bd{\mu}_n,\bd{s}_n \rangle^{\bullet}}, 
\end{equation}
\begin{equation}
\dot{t}_i.|\bd{\nu}_n,\bd{s}_n \rangle^{\bullet} = p^{M_i(\bd{\nu}_n;\bd{s}_n;l)}\, |\bd{\nu}_n,\bd{s}_n \rangle^{\bullet},
\end{equation}
\begin{equation}
\dot{\partial}.|\bd{\nu}_n,\bd{s}_n \rangle^{\bullet} = -\bigl( \Delta(\bd{s}_n,l) + N_0(\bd{\nu}_n;\bd{s}_n;l) \bigr) 
\,|\bd{\nu}_n,\bd{s}_n \rangle^{\bullet}.
\end{equation} \cqfd
\end{thm}

Theorems \ref{thm_action_Uq_sur_Lambda_s} and \ref{thm_action_Up_sur_Lambda_s} show in particular that the vectors of the standard basis
of $\Lambda^s$ are weight vectors for the actions of $\Uq$ and $\Up$, and the weights are given by:
\begin{cor}[\cite{U}] \label{cor_formulas_weights} With obvious notation, we have 
\begin{eqnarray}
\hspace{1cm} \wt(|\bd{\lambda}_l, \bd{s}_l \rangle) &=&
 - \Delta(\bd{s}_l,n) \delta + \Lambda_{s_1}+ \cdots + \Lambda_{s_l} - \ds \sum_{i=0}^{n-1} N_i(\bd{\lambda}_l;\bd{s}_l;n) \,\alpha_i, 
\label{poidsl} \\
\hspace{1cm} \dot{\wt}(|\bd{\lambda}_l, \bd{s}_l \rangle) &=& 
 - \bigl(\Delta(\bd{s}_l,n)+ N_0(\bd{\lambda}_l;\bd{s}_l;n) \bigr) \dot{\delta} 
 + (n-s_1+s_l) \dot{\Lambda}_0+ \ds \sum_{i=1}^{l-1}{(s_i-s_{i+1})  \,\dot{\Lambda}_i},
\label{poidsptl} \\
\hspace{1cm} \dot{\wt}(|\bd{\lambda}_n, \bd{s}_n \rangle^{\bullet}) &=& 
- \Delta(\bd{s}_n,l) \dot{\delta} + \dot{\Lambda}_{s_1}+ \cdots + \dot{\Lambda}_{s_n} - \ds \sum_{i=0}^{l-1} N_i(\bd{\lambda}_n;\bd{s}_n;l) 
\,\dot{\alpha}_i, 
\label{poidsptn} \\
\hspace{1cm} \wt (|\bd{\lambda}_n, \bd{s}_n \rangle^{\bullet}) &=& 
 - \bigl(\Delta(\bd{s}_n,l)+ N_0(\bd{\lambda}_n;\bd{s}_n;l) \bigr) \delta 
  + (l-s_1+s_n) \Lambda_0 + \ds \sum_{i=1}^{n-1}{(s_i-s_{i+1}) \,\Lambda_i}. 
\label{poidsn}
\end{eqnarray} \cqfd
\end{cor}

\begin{definition}
For $m \in \mZ^*$, define an endomorphism $B_m$ of $\Lambda^s$ by

\begin{equation} 
B_m(u_{k_1} \wedge u_{k_2} \wedge \cdots) := \ds \sum_{j=1}^{+\infty} { u_{k_1} \wedge \cdots \wedge u_{k_{j-1}} \wedge u_{k_j-nlm} \wedge 
u_{k_{j+1}} \wedge \cdots} \qquad \bigl( (k_1,k_2,\ldots) \in P^{++}(s) \bigr).
\end{equation} \fini
\end{definition}
Note that by \cite[Lemma 3.18]{U}, the sum above involves only finitely many nonzero terms, hence $B_m$ is well-defined.
This definition comes from a passage to the limit $r \rightarrow \infty$ in the action of the center of the Hecke algebra of 
$\widehat{\mathfrak{S}}_r$ on $q$-wedge products of $r$ factors. However, the operators $B_m$ do not commute, but by 
\cite[Prop. 4.4 \& 4.5]{U}, they span a Heisenberg algebra
\begin{equation}
\mathcal{H} := \langle B_m \mid m \in \mZ^* \rangle.
\end{equation}

We now recall some results concerning the actions of $\Uq$, $\Up$ and $\mathcal{H}$ on $\Lambda^s$. 

\begin{prop}[\cite{U}] The actions of $\Uqprime$, $\Upprime$ and $\mathcal{H}$ on $\Lambda^s$ pairwise commute.\hspace{-1em} \cqfd
\end{prop} 

For $L$, $N \in \mN^*$, introduce the finite set
\begin{equation} \label{definition_A_l_n_s}
A_{L,N}(s):= \set{(r_1,\ldots,r_L) \in \mZ^L(s) \mid r_1 \geq \cdots \geq r_L, \ r_1-r_L \leq N}.
\end{equation}
It is not hard to see that if $\bd{r}_l \in \mZ^l(s)$ and $\bd{r}_n \in \mZ^n(s)$ are such that 
$|\bd{\emptyset}_l,\bd{r}_l \rangle = |\bd{\emptyset}_n,\bd{r}_n \rangle^{\bullet}$, then $\bd{r}_l \in A_{l,n}(s)$ and 
$\bd{r}_n \in A_{n,l}(s)$. Conversely, if $\bd{r}_l \in A_{l,n}(s)$, then there exists a unique $\bd{r}_n \in A_{n,l}(s)$ such that
$|\bd{\emptyset}_n,\bd{r}_n \rangle^{\bullet} = |\bd{\emptyset}_l,\bd{r}_l \rangle$, and if $\bd{r}_n \in A_{n,l}(s)$, then there exists a unique 
$\bd{r}_l \in A_{l,n}(s)$ such that $|\bd{\emptyset}_l,\bd{r}_l \rangle = |\bd{\emptyset}_n,\bd{r}_n \rangle^{\bullet}$. Therefore, 
$$\set{|\bd{\emptyset}_l,\bd{r}_l \rangle \mid \bd{r}_l \in A_{l,n}(s)}
=\set{|\bd{\emptyset}_n,\bd{r}_n \rangle^{\bullet} \mid \bd{r}_n \in A_{n,l}(s)}$$ is a set of highest weight vectors simultaneously 
for the actions of $\Uqprime$ and $\Upprime$. It is easy to see that these vectors are also singular for the action of $\mathcal{H}$, that is, 
they are annihilated by the $B_m$, $m>0$. It turns out that these vectors are the only singular vectors simultaneously for the actions of 
$\Uqprime$, $\Upprime$ and $\mathcal{H}$, and we have the following theorem.      

\begin{thm}[\cite{U}, Thm. 4.8] \label{thm_dec_Lambda_s} We have
$$
\Lambda^s = \ds \bigoplus_{\bd{r}_l \in A_{l,n}(s)}{\Uqprime \otimes \mathcal{H} \otimes \Upprime.  
|\bd{\emptyset}_l,\bd{r}_l} \rangle
= \ds \bigoplus_{\bd{r}_n \in A_{n,l}(s)}{\Uqprime \otimes \mathcal{H} \otimes \Upprime.  
|\bd{\emptyset}_n,\bd{r}_n \rangle^{\bullet}}.$$ \cqfd
\end{thm}

\subsubsection{The involution $\barre$ of $\Lambda^s$}

Following \cite{U}, the space $\Lambda^s$ can be endowed with an involution $\barre$. Instead of recalling the definition of this 
involution, we give its main properties (they turn out to characterize it completely).

\begin{prop}[\cite{U}] \label{prop_def_barre} There exists an involution $\barre$ of $\Lambda^s$ such that:
\begin{itemize}
\item[\rm (i)] $\barre$ is a $\mQ$-linear map of $\Lambda^s$ such that for all $u \in \Lambda^s$, $k \in \mZ$, we have $\overline{q^k\, u}=q^{-k}\,\overline{u}$.
\item[\rm (ii)] \emph{(Unitriangularity property)}. For all $\lambda \in \Pi$, we have 
$$\overline{| \lambda,s \rangle} \in | \lambda,s \rangle + \ds \bigoplus_{\mu \lhd \lambda}{\mZ[q,q^{-1}] \,|\mu,s \rangle},$$ where $\lhd$ stands 
for the dominance ordering on partitions.
\item[\rm (iii)] For all $\lambda \in \Pi$, we have $\wt(\overline{| \lambda,s \rangle})=\wt(|\lambda,s \rangle)$ and 
$\wtpt(\overline{| \lambda,s \rangle})=\wtpt(|\lambda,s \rangle)$.
\item[\rm (iv)] For all $0 \leq i \leq n-1$,  $0 \leq j \leq l-1$, $m<0$, $v \in \Lambda^s$, we have
$$\overline{f_i.v}=f_i.\overline{v}, \quad \overline{\dot{f}_j.v}=\dot{f}_j.\overline{v} \quad \mbox{and} \quad
\overline{B_m.v}=B_m.\overline{v}.$$
\end{itemize} 
\end{prop}

\begin{proof} Let $\barre$ be the involution of $\Lambda^s$ defined in \cite[Prop. 3.23 \& Eq. (39)]{U}. By construction, (i) holds. The other
statements come from \cite[Prop. 4.11 \& 4.12]{U} and Corollary \ref{cor_formulas_weights}.
\end{proof}

\subsection{$q$-deformed higher-level Fock spaces}
\subsubsection{Definition}
By Theorem \ref{thm_action_Uq_sur_Lambda_s}, the space
\begin{equation}
\Fq:= \ds \bigoplus_{\bd{\lambda}_l \in \Pi^l}{\mQ(q) |\bd{\lambda}_l, \bd{s}_l \rangle} \subset \Lambda^s \qquad (\bd{s}_l \in \mZ^l(s))
\end{equation}
is a $\Uq$-submodule of $\Lambda^s$. The reader should be aware that $\Fq$ is \emph{not} a $\Up$-submodule of $\Lambda^s$.  
In a similar way, by Theorem \ref{thm_action_Up_sur_Lambda_s}, the space
\begin{equation}
\Fp:= \ds \bigoplus_{\bd{\lambda}_n \in \Pi^n}{\mQ(q) |\bd{\lambda}_n, \bd{s}_n \rangle}^{\bullet} \subset \Lambda^s \qquad (\bd{s}_n \in \mZ^n(s))
\end{equation}
is a $\Up$-submodule of $\Lambda^s$. 

\begin{definition}[\cite{U}]
The representations $\Fq$ and $\Fp$ ($\bd{s}_l \in \mZ^l(s)$, $\bd{s}_n \in \mZ^n(s)$) are called \emph{($q$-deformed) Fock spaces}. When $l>1$ and 
$n>1$, we speak of higher-level Fock spaces. \fini
\end{definition}

Since the maps $\tau_l$ and $\tau'_n$ are bijections, we have the following decompositions:

\begin{equation}
\Lambda^s = \ds \bigoplus_{\bd{s}_l \in \mZ^l(s)} \Fq = \ds \bigoplus_{\bd{s}_n \in \mZ^n(s)} \Fp.
\end{equation} 
Neither of these decompositions is compatible with the decompositions of $\Lambda^s$ given in Theorem \ref{thm_dec_Lambda_s}.

\subsubsection{Fock spaces as weight subspaces of $\Lambda^s$. Actions of the Weyl groups.}

Let $N$, $L \in \mN^*$. Recall the definition of $\mQ^{L}(s)$ and $\mQ^{L}(N)$ from (\ref{eq_Z_l_s}) and define a map
\begin{equation}
\theta_{L,N} \, : \mQ^{L}(s) \rightarrow \mQ^{L}(N), \quad (s_1,\ldots,s_L) \mapsto (N-s_1+s_L,s_1-s_2,\ldots, s_{L-1}-s_L).
\end{equation}
It is easy to see that $\theta_{L,N}$ is bijective. Moreover, for $(a_1,\ldots,a_L) \in \mQ^{L}(N)$, the $l$-tuple
$(s_1,\ldots,s_L):=\theta_{L,N}^{-1}(a_1,\ldots,a_L)$ is given by
\begin{equation}
s_i = \ds \frac{1}{L} \left( s - \sum_{j=1}^{L-1} ja_j \right) + \sum_{j=i+1}^{L} a_j \qquad (1 \leq i \leq L).
\end{equation}
The next result shows that the Fock spaces are sums of certain weight subspaces of $\Lambda^s$.
The proof follows easily from Corollary \ref{cor_formulas_weights}.
\newpage 
\begin{prop}[\cite{U}] \ \\[-5mm] \label{fock_spaces_weight_subspaces}
\begin{itemize}
\item[\rm (i)] Let $\bd{s}_n \in \mZ^n(s)$. Let $(a_0,\ldots,a_{n-1}):=\theta_{n,l}(\bd{s}_n)$ and $w := \sum_{i=0}^{n-1}{a_i \Lambda}_i$.
Then $$\Fp= \ds \bigoplus_{d \in \mZ}{\Lambda^s \langle w+d \delta} \rangle.$$
\item[\rm (ii)] Let $\bd{s}_l \in \mZ^l(s)$. Let $(a_0,\ldots,a_{l-1}):=\theta_{l,n}(\bd{s}_l)$ and $\dot{w}:=\sum_{i=0}^{l-1}{a_i \dot{\Lambda}_i}$.
Then $$\Fq= \ds \bigoplus_{d \in \mZ}{\Lambda^s \langle {\dot{w}+d \dot{\delta}}} \rangle.$$ \cqfd
\end{itemize}
\end{prop}

Note that the operator $B_m$ $(m \in \mZ^*)$ maps the weight subspace $\Lambda^s \langle w \rangle$ (\resp $\Lambda^s \langle \dot{w} \rangle$) into 
$\Lambda^s \langle w+m \delta \rangle$ (\resp $\Lambda^s \langle \dot{w}+m \dot{\delta} \rangle$). Therefore, by Proposition
\ref{fock_spaces_weight_subspaces}, the Fock spaces $\Fq$ and $\Fp$ ($\bd{s}_l \in \mZ^l(s)$, $\bd{s}_n \in \mZ^n(s)$) are stable under the 
action of $\mathcal{H}$. \\ 

We now compare some weight subspaces of the Fock spaces. The proof follows again from Corollary \ref{cor_formulas_weights}. 

\begin{prop} \label{prop_bij_pds_charges}
Let $\bd{s}_l \in \mZ^l(s)$ and $w$ be a weight of $\Fq$. Then there exists a unique pair $(\bd{s}_n,\dot{w})$ such that
$\bd{s}_n \in \mZ^n(s)$, $\dot{w}$ is a weight of $\Fp$ and $\Fq \langle w \rangle = \Fp \langle \dot{w} \rangle$. More precisely, write
$w = d \delta + \sum_{i=0}^{n-1}{a_i \Lambda_i}$ with $a_0,\ldots,a_{n-1},d \in \mZ$, 
$\bd{s}_l $, and put $s_0:=n+s_l$. Then we have 
$\bd{s}_n=\theta_{n,l}^{-1}(a_0,\ldots,a_{n-1})$ and $\dot{w}= d \dot{\delta} + \sum_{i=0}^{l-1}{(s_i-s_{i+1}) \dot{\Lambda}_i}.$ \cqfd
\end{prop}

\begin{example} Take $n=3$, $l=2$, $\bd{s}_l=(1,0)$ and $w = -2 \Lambda_0 + \Lambda_1 + 3 \Lambda_2 - 2 \delta$. Then by (\ref{poidsl}),
we have $\wt \big(\big|\bigl((1,1),(1)\bigr),\bd{s}_l \big\rangle \big)=w$, so $w$ is a weight of $\Fq$.
By Proposition \ref{prop_bij_pds_charges}, we have $\Fq \langle w \rangle = \Fp \langle \dot{w} \rangle$ with $\bd{s}_n=(2,1,-2)$ and 
$\dot{w} = 2 \dot{\Lambda}_0+\dot{\Lambda}_1-2 \dot{\delta}$. Moreover, using
(\ref{poidsl}) and (\ref{poidsptn}), we see that for all 
$|\bd{\lambda}_l,\bd{s}_l \rangle = |\bd{\lambda}_n,\bd{s}_n \rangle^{\bullet} \in \Fq \langle w \rangle = \Fp \langle \dot{w} \rangle$,
we have $N_0(\bd{\lambda}_l;\bd{s}_l;n)=2$, $N_1(\bd{\lambda}_l;\bd{s}_l;n)=1$, $N_2(\bd{\lambda}_l;\bd{s}_l;n)=0$ and
$N_0(\bd{\lambda}_n;\bd{s}_n;l)=N_1(\bd{\lambda}_n;\bd{s}_n;l)=0$ (this shows \emph{a posteriori} that 
$\operatorname{dim} \Fq \langle w \rangle =1$). \fini
\end{example}

We now deal with the actions of the Weyl groups of $\Uq$ and $\Up$ on the set of the weight subspaces of $\Lambda^s$. Recall that 
$W_n = \langle \sigma_0,\ldots,\sigma_{n-1} \rangle$ is the Weyl group of $\Uq$. By (\ref{def_Weyl_group}), $W_n$ acts on the
weight lattice $\bigoplus_{i=0}^{n-1}{\mZ \Lambda_i} \oplus \mZ \delta$ by
\begin{equation} \label{eq_action_Weyl_group_fundamental_weights}
\sigma_i.\delta = \delta \quad \mbox{and} \quad
 \sigma_i.\Lambda_j= \left\{
 \begin{array}{lc}
 \Lambda_j & \mbox{ if } j \neq i, \\
 \Lambda_{i-1}+ \Lambda_{i+1} - \Lambda_i -\delta_{i,0} \, \delta & \mbox{ if } j = i \\
\end{array} \right. \qquad (0 \leq i,j \leq n-1).
\end{equation}
 
Moreover, it is easy to see that $W_n$ acts faithfully on $\mZ^n(s)$ by
\begin{equation}
\left\{ \begin{array}{rcl}
\sigma_0.(s_1,\ldots,s_n) &=& (s_n+l,s_2,\ldots,s_{n-1},s_1-l), \\
\sigma_i.(s_1,\ldots,s_n) &=& (s_1,\ldots,s_{i+1},s_i,\ldots,s_n) \qquad (1 \leq i \leq n-1),\\
\end{array} \right.
\end{equation}
and the set $A_{n,l}(s)$ defined by (\ref{definition_A_l_n_s}) is a fundamental domain for this action. In a similar way, one can define two actions 
of the Weyl group $\dot{W}_l$ of $\Up$, one on the weight lattice $\bigoplus_{i=0}^{l-1}{\mZ \dot{\Lambda}_i} \oplus \mZ \dot{\delta}$ and one on 
$\mZ^l(s)$. The following lemma will be useful later.

\begin{lemma} \label{lemma_action_gpe_Weyl_on_Fock_spaces}
Let $\bd{s}_n \in \mZ^n(s)$ and $\dot{w}$ be a weight of $\Fp$. Let $(\bd{s}_l,w) \in \mZ^l(s) \times \mathcal{P}(\Lambda^s)$ be the
unique pair such that $\Fp \langle \dot{w} \rangle = \Fq \langle w \rangle$ (see Proposition \ref{prop_bij_pds_charges}). Let 
$\dot{\sigma} \in \dot{W}_l$. In the same way, let $(\bd{t}_l,w') \in \mZ^l(s) \times \mathcal{P}(\Lambda^s)$ be the
unique pair such that $\Fp \langle \dot{\sigma}.\dot{w} \rangle = \Fq[\bd{t}_l] \langle w' \rangle$. Then we have
$$\bd{t}_l=\dot{\sigma}.\bd{s}_l \qquad \mbox{and} \qquad w'=w+\wt(|\bd{\emptyset}_l, \bd{t}_l \rangle)-\wt(|\bd{\emptyset}_l, \bd{s}_l \rangle).$$
\end{lemma}

\begin{proof} The proof follows immediately from the formulas given in Proposition \ref{prop_bij_pds_charges}. \end{proof}

\subsubsection{The lower crystal basis $(\mathcal{L}[\bd{s}_l],\mathcal{B}[\bd{s}_l])$ of $\Fq$ at $q=0$}

Let $\bd{s}_l \in \mZ^l(s)$. Following \cite{Kas1}, let $\mA \subset \mQ(q)$ be the ring of rational functions which are regular at $q=0$, 
$\mathcal{L}[\bd{s}_l]:=\bigoplus_{\bd{\lambda}_l \in \Pi^l}{\mA \, |\bd{\lambda}_l,\bd{s}_l \rangle}$ and
for $0 \leq i \leq n-1$, let $\eilow$, $\filow$, $\eiup$ and $\fiup$ 
denote Kashiwara's operators acting on $\mathcal{L}[\bd{s}_l]$. The following lemma shows that sometimes, certain powers of the operators 
$\eilow$ and $\eiup$ coincide (one has an analogous result for $\filow$ and $\fiup$).

\begin{lemma} \label{lemma_Kashiwara_operators}
Let $\bd{s}_l \in \mZ^l(s)$, $w \in \mathcal{P}(\Fq)$, $u \in (\operatorname{Ker} e_i) \cap \Fq \langle w \rangle$ and $k:=(w,\alpha_i)$.
Then we have $$(\eiup)^k.(f_i^{(k)}.u)=(\eilow)^k.(f_i^{(k)}.u)=u.$$
\end{lemma} 

\begin{proof} 
The second equality follows easily by induction on $k$ from the definition of $\eilow$. Let us now show that
$(\eiup)^k.(f_i^{(k)}.u)=u$. Note that for $0 \leq j \leq k$, we have 
$$(\wt(f_i^{(k-j)}.u),\alpha_i) = (\wt(u),\alpha_i) - (k-j)(\alpha_i,\alpha_i) = 2j-k.$$ 
By induction on $0 \leq k' \leq k$, we get therefore, by definition of $\eiup$,
$$(\eiup)^{k'}.(f_i^{(k)}.u) = \ds \left( {\prod_{j=0}^{k'-1} \frac{[(2j-k)+(k-j)+1]}{[k-j]}} \right) \, f_i^{(k-k')}.u \qquad (0 \leq k' \leq k).$$
As a consequence, we have $(\eiup)^{k}.(f_i^{(k)}.u)= \ds \left( {\prod_{j=0}^{k-1} \frac{[j+1]}{[k-j]}} \right) \, u = u$.
\end{proof} 

Put
\begin{equation}
\mathcal{B}[\bd{s}_l]:=\set{|\bd{\lambda}_l,\bd{s}_l \rangle \bmod q \mathcal{L}[\bd{s}_l] \mid \bd{\lambda}_l \in \Pi^l}
\end{equation} 
In the sequel, we shall write more briefly $\bd{\lambda}_l$ for the
element in $\mathcal{B}[\bd{s}_l]$ indexed by the corresponding multi-partition. By \cite{JMMO},\cite{FLOTW},\cite{U}, the pair 
$(\mathcal{L}[\bd{s}_l],\mathcal{B}[\bd{s}_l])$ is a lower crystal basis of $\Fq$ at $q=0$ in the sense of \cite{Kas1}, and the
crystal graph $\mathcal{B}[\bd{s}_l]$ 
contains the arrow $\bd{\lambda}_l \stackrel{i}{\longrightarrow} \bd{\mu}_l$ if and only if the multi-partition $\bd{\mu}_l$ is obtained 
from $\bd{\lambda}_l$ by adding a good $i$-node in the sense of \cite[Thm. 2.4]{U}. We shall still denote by $\eilow$ and $\filow$ Kashiwara's
operators acting on $\mathcal{B}[\bd{s}_l] \cup \set{0}$.   

\subsubsection{Uglov's canonical bases of the Fock spaces}

By Propositions \ref{prop_def_barre} (iii) and \ref{fock_spaces_weight_subspaces}, the Fock spaces $\Fq$ and $\Fp$ 
($\bd{s}_l \in \mZ^l(s)$, $\bd{s}_n \in \mZ^n(s)$) are stable under the involution $\barre$. The involution induced on these spaces will still be 
denoted by $\barre$. Let $\bd{s}_l \in \mZ^l(s)$. For $\bd{\mu}_l \in \Pi^l$, write 
\begin{equation}
\overline{|\bd{\mu}_l,\bd{s}_l \rangle} = \ds \sum_{\bd{\lambda}_l \in \Pi^l}
{a_{\bd{\lambda}_l,\bd{\mu}_l ;\, \bd{s}_l}(q)\, |\bd{\lambda}_l,\bd{s}_l \rangle}
\end{equation}
\noindent with $a_{\bd{\lambda}_l,\bd{\mu}_l ;\,\bd{s}_l}(q) \in \mZ[q,q^{-1}]$, and let 
\begin{equation}
A_{\bd{s}_l}(q):=\bigl(a_{\bd{\lambda}_l,\bd{\mu}_l ;\,\bd{s}_l}(q) \bigr)_{\bd{\lambda}_l,\bd{\mu}_l \in \Pi^l}
\end{equation}
denote the matrix of the involution $\barre$ of $\Fq$. Since the weight subspaces of $\Fq$ are stable under the involution $\barre$, 
(\ref{poidsl}) implies that $a_{\bd{\lambda}_l,\bd{\mu}_l ;\,\bd{s}_l}(q)$ is zero unless $|\bd{\lambda}_l| = |\bd{\mu}_l|$,
where $|\bd{\lambda}_l|$ (\resp $|\bd{\mu}_l|$) denotes the number of boxes contained in the Young diagram of $|\bd{\lambda}_l|$ 
(\resp $|\bd{\mu}_l|$). By Proposition \ref{prop_def_barre} (ii), the matrix $A_{\bd{s}_l}(q)$ is unitriangular. One can therefore define, by a 
classical argument, canonical bases of the Fock space $\Fq$ as follows. 

\begin{thm}[\cite{U}] \label{thm_bases_canoniques_Fqsl}
Let $\bd{s}_l \in \mZ^l(s)$. Then there exists a unique base 
$$\set{G^+(\bd{\lambda}_l, \bd{s}_l) \mid \bd{\lambda}_l \in \Pi^l} \qquad \Bigl( \mbox{\resp}
\set{ G^-(\bd{\lambda}_l, \bd{s}_l) \mid \bd{\lambda}_l \in \Pi^l} \Bigr)$$
of $\Fq$ such that:
$$\begin{array}{lll}
\mbox{\rm (i)} & \overline{G^+(\bd{\lambda}_l, \bd{s}_l)}=G^+(\bd{\lambda}_l, \bd{s}_l) & (\mbox{\resp} 
 \overline{G^-(\bd{\lambda}_l, \bd{s}_l)}=G^-(\bd{\lambda}_l, \bd{s}_l) \; \mbox{\emph{),}}  \\ 
\mbox{\rm (ii)} & G^+(\bd{\lambda}_l, \bd{s}_l) \equiv |\bd{\lambda}_l, \bd{s}_l \rangle \bmod \ q\mathcal{L}^{+}[\bd{s}_l] & (\mbox{\resp}
G^-(\bd{\lambda}_l, \bd{s}_l) \equiv |\bd{\lambda}_l, \bd{s}_l \rangle \bmod \ q^{-1}\mathcal{L}^{-}[\bd{s}_l] \; ), \\ 
\end{array}$$
where $\mathcal{L}^{\epsilon}[\bd{s}_l]:= \ds \bigoplus_{\bd{\lambda}_l \in \Pi^l}{\mZ[q^\epsilon]\,|\bd{\lambda}_l, \bd{s}_l \rangle}$ 
$(\epsilon=\pm 1)$. \cqfd
\end{thm}

\begin{definition}
The bases $\set{G^+(\bd{\lambda}_l, \bd{s}_l) \mid \bd{\lambda}_l \in \Pi^l}$ and $\set{G^-(\bd{\lambda}_l,\bd{s}_l) \mid \bd{\lambda}_l \in \Pi^l}$
are called \emph{Uglov's canonical bases} of $\Fq$. Define entries $\Delta^{+}_{\bd{\lambda}_l,\bd{\mu}_l ;\,\bd{s}_l}(q)$,
$\Delta^{-}_{\bd{\lambda}_l,\bd{\mu}_l ;\,\bd{s}_l}(q) \in \mZ[q,q^{-1}]$ ($\bd{\lambda}_l$, $\bd{\mu}_l \in \Pi^l$) by
\begin{equation} \label{def_Delta_q}
G^{+}(\bd{\mu}_l,\bd{s}_l) = \ds \sum_{\bd{\lambda}_l \in \Pi^l}
{\Delta^{+}_{\bd{\lambda}_l,\bd{\mu}_l ;\,\bd{s}_l}(q) \,|\bd{\lambda}_l, \bd{s}_l \rangle}, \quad
G^{-}(\bd{\mu}_l,\bd{s}_l) = \ds \sum_{\bd{\lambda}_l \in \Pi^l}
{\Delta^{-}_{\bd{\lambda}_l,\bd{\mu}_l ;\,\bd{s}_l}(q) \,|\bd{\lambda}_l, \bd{s}_l \rangle},
\end{equation}
and denote by 
\begin{equation}
\Delta^{\epsilon}_{\bd{s}_l}(q):=\bigl(\Delta^{\epsilon}_{\bd{\lambda}_l,\bd{\mu}_l ;\,\bd{s}_l}(q) \bigr)_{\bd{\lambda}_l,\bd{\mu}_l \in \Pi^l}
\qquad  (\epsilon = \pm 1)
\end{equation}
the transition matrices between the standard and the canonical bases of $\Fq$. \fini
\end{definition}

By \cite{U}, the entries of $\Delta^{+}_{\bd{s}_l}(q)$ (\resp $\Delta^{-}_{\bd{s}_l}(q)$) are Kazhdan-Lusztig polynomials of parabolic
submodules of affine Hecke algebras of type $A$, so by \cite{KT}, these polynomials are in $\mN[q]$ (\resp $\mN[p]$). Moreover, both canonical 
bases of $\Fq$ are dual to each other with respect to a certain bilinear form, which gives an inversion formula for Kazhdan-Lusztig polynomials; see 
\cite[Thm. 5.15]{U}. By \cite{U}, the basis 
$\set{G^+(\bd{\lambda}_l, \bd{s}_l) \mid \bd{\lambda}_l \in \Pi^l}$ is a lower global crystal basis (in the sense of \cite{Kas1}) of the integrable 
$\Uq$-module $\Fq$. \\

Let $\bd{s}_n \in \mZ^n(s)$. In a similar way, one can define canonical bases 
$\set{G^{\epsilon}(\bd{\lambda}_n, \bd{s}_n)^{\bullet} \mid \bd{\lambda}_n \in \Pi^n}$
$(\epsilon = \pm 1)$ of the Fock space $\Fp$. By \cite{U}, the basis 
$\set{G^{-}(\bd{\lambda}_n, \bd{s}_n)^{\bullet} \mid \bd{\lambda}_n \in \Pi^n}$ is a lower global crystal basis of the integrable 
$\Up$-module $\Fp$. For $\bd{\mu}_n \in \Pi^n$, $\epsilon = \pm 1$, write
\begin{equation}
G^{\epsilon}(\bd{\mu}_n,\bd{s}_n)^{\bullet} = \ds \sum_{\bd{\lambda}_n \in \Pi^n}
{\dot{\Delta}^{\epsilon}_{\bd{\lambda}_n,\bd{\mu}_n;\, \bd{s}_n}(q)\,|\bd{\lambda}_n, \bd{s}_n \rangle}^{\bullet},
\end{equation}
where the entries $\dot{\Delta}^{\epsilon}_{\bd{\lambda}_n,\bd{\mu}_n;\, \bd{s}_n}(q)$ are in $\mZ[q,q^{-1}]$. Since the weight subspaces of 
$\Lambda^s$ are stable under the involution $\barre$, we have $\dot{\Delta}^{\epsilon}_{\bd{\lambda}_n,\bd{\mu}_n ;\, \bd{s}_n}(q) = 0$ unless 
$\wt(|\bd{\lambda}_n,\bd{s}_n \rangle^{\bullet})=\wt(|\bd{\mu}_n,\bd{s}_n \rangle^{\bullet})$ and
$\wtpt(|\bd{\lambda}_n,\bd{s}_n \rangle^{\bullet})=\wtpt(|\bd{\mu}_n,\bd{s}_n \rangle^{\bullet})$. In this case, by
Corollary \ref{cor_formulas_weights}, there exist $\bd{s}_l \in \mZ^l(s)$ and $\bd{\lambda}_l$, $\bd{\mu}_l \in \Pi^l$ such that
$|\bd{\lambda}_l,\bd{s}_l \rangle = |\bd{\lambda}_n,\bd{s}_n \rangle^{\bullet}$ and 
$|\bd{\mu}_l,\bd{s}_l \rangle = |\bd{\mu}_n,\bd{s}_n \rangle^{\bullet}$. It is then not hard to see that
\begin{equation} \label{equation_lien_Delta_Deltapoint} 
\dot{\Delta}^{\epsilon}_{\bd{\lambda}_n,\bd{\mu}_n ;\, \bd{s}_n}(q) = \Delta^{\epsilon}_{\bd{\lambda}_l,\bd{\mu}_l ;\,\bd{s}_l}(q).
\end{equation}

\section{Comparison of canonical bases of weight subspaces of $\Fq$ $(\bd{s_l} \in \mZ^l(s))$}

In the sequel we use the following notation.
\begin{notation} For $\bd{s}_l \in \mZ^l(s)$ and $w \in \mathcal{P}(\Fq)$, put
\begin{equation}
\Pi^l(\bd{s}_l;w):= \set{ \bd{\lambda}_l \in \Pi^l \ \big| \ |\bd{\lambda}_l,\bd{s}_l \rangle \in \Lambda^s \langle w \rangle},
\end{equation}
and define similarly $\Pi^n(\bd{s}_n;\dot{w})$ for $\bd{s}_n \in \mZ^n(s)$ and $\dot{w} \in \dot{\mathcal{P}}(\Fp)$. \fini
\end{notation}

\begin{definition} \label{def_bc_semblables}
Let $\bd{s}_l$, $\bd{t}_l \in \mZ^l(s)$, $w \in \mathcal{P}(\Fq)$ and $w' \in \mathcal{P}(\Fq[\bd{t}_l])$. 
We say that the canonical bases of $\Fq \langle w \rangle$ and $\Fq[\bd{t}_l] \langle w' \rangle$ are \emph{similar} if there exists a bijection 
$$\sigma : \Pi^l(\bd{s}_l;w) \rightarrow \Pi^l(\bd{t}_l;w')$$ such that for all $\bd{\lambda}_l$, $\bd{\mu}_l \in \Pi^l(\bd{s}_l;w)$,
$\epsilon = \pm 1$, we have 
$$\Delta^{\epsilon}_{\sigma(\bd{\lambda}_l),\sigma(\bd{\mu}_l) ;\, \bd{s}_l}(q)=\Delta^{\epsilon}_{\bd{\lambda}_l,\bd{\mu}_l ;\,\bd{s}_l}(q).$$
In other words, the canonical bases of $\Fq \langle w \rangle$ and $\Fq[\bd{t}_l] \langle w' \rangle$ are similar if the transition matrices 
between the standard bases and the canonical bases are equal up to a reindexing of rows and columns. \fini
\end{definition}

\begin{notation}
Throughout this section we fix a multi-charge $\bd{s_l} \in \mZ^l(s)$ and an integer $i \in \interv{0}{n-1}$. To simplify, we drop the 
multi-charge $\bd{s}_l$ in the notation of this section, that is we denote by $\bd{\lambda}_l$ (\resp $G^{\pm}(\bd{\lambda}_l)$)
the vector of the standard (\resp canonical) basis of $\Fq$ indexed by the corresponding multi-partition and so on. In particular, we use the
notation $\bd{\lambda}_l$ either for a vector of the standard basis of $\Fq$ or for a vertex in the crystal graph 
$\mathcal{B}:=\mathcal{B}[\bd{s}_l]$. \fini
\end{notation}

\subsection{The bijection $\sigma_i$} \label{paragraphe_bij_scopes}
We recall here the definition of the involution $\sigma_i$ of the crystal graph $\mathcal{B}$. We sometimes view $\sigma_i$ as a bijection of
$\Pi^l$.

\begin{definition} \label{def_sigma_i} 
Let $\bd{\lambda}_l \in \mathcal{B} \cong \Pi^l$ and $i \in \interv{0}{n-1}$. Let $\mathcal{C}$ be the $i$-chain in $\mathcal{B}$ containing 
$\bd{\lambda}_l$. Let $\sigma_i(\bd{\lambda}_l) \in \mathcal{B} \cong \Pi^l$ be the unique element in $\mathcal{C}$ such that
$\wt(\sigma_i(\bd{\lambda}_l))=\sigma_i.(\wt(\bd{\lambda}_l))$. In other words, $\sigma_i(\bd{\lambda}_l)$ is obtained from $\bd{\lambda}_l$ via a 
central symmetry in the middle of $\mathcal{C}$. This defines an involution $\sigma_i$ of $\mathcal{B}$. This map induces, for 
$w \in \mathcal{P}(\Fq)$, a bijection 
\begin{equation} \label{def_sigma_i_2}
\sigma_i : \Pi^l(\bd{s}_l;w) \stackrel{\sim}{\longrightarrow} \Pi^l(\bd{s}_l;\sigma_i.w).
\end{equation} \fini 
\end{definition}

By \cite{Kas2}, the definition of $\sigma_0,\ldots,\sigma_{n-1}$ as bijections of 
$\mathcal{B}$ gives actually rise to an action of the Weyl group $W_n$ on $\mathcal{B}$, but we do not need this fact in the sequel. 

\begin{prop} \label{bij_scopes} Let $w \in \mathcal{P}(\Fq)$ and $i \in \interv{0}{n-1}$ be such that $w + \alpha_i \notin \mathcal{P}(\Fq)$. 
Let $\bd{\lambda}_l \in \Pi^l(\bd{s}_l;w)$ and $\bd{\mu}_l:=\sigma_i(\bd{\lambda}_l)$. Then we have the following:

\begin{itemize}
\item[\rm (i)] $\bd{\mu}_l$ is the multi-partition obtained by adding to $\bd{\lambda}_l$ all its addable $i$-nodes, and we have 
$|\bd{\mu}_l \setminus \bd{\lambda}_l|=k_i:= (w,\alpha_i)$.
\item[\rm (ii)] In $\Fq$, we have $\bd{\mu}_l = f_i^{(k_i)}. \bd{\lambda}_l$ and $\bd{\lambda}_l = e_i^{(k_i)}. \bd{\mu}_l$.
\item[\rm (iii)] In $\mathcal{B}$, we have $\bd{\mu}_l = (\filow)^{k_i}.\bd{\lambda}_l$ and $\bd{\lambda}_l = (\eilow)^{k_i}.\bd{\mu}_l$. 
\end{itemize}
\end{prop}

\begin{proof} Let $\mathcal{C}$ be the $i$-chain in $\mathcal{B}$ containing $\bd{\lambda}_l$ and $\bd{\mu}_l$. Since
$w + \alpha_i$ is not a weight of $\Fq$, $\bd{\lambda}_l$ has no removable $i$-node. This implies, by \cite[Thm. 2.4]{U},
that $\bd{\lambda}_l$ is the head of the chain $\mathcal{C}$, and by symmetry $\bd{\mu}_l$ is the tail of $\mathcal{C}$.
Note that since $w + \alpha_i$ is not a weight of $\Fq$, $\sigma_i.(w + \alpha_i)=(\sigma_i.w)-\alpha_i$ is also not a weight of $\Fq$,
so $\bd{\mu}_l$ has no addable $i$-node. By \cite[Thm. 2.4]{U}, $\bd{\mu}_l$ is obtained by adding some $i$-nodes (let us say $k$ of them) to 
$\bd{\lambda}_l$. The integer $k=|\bd{\mu}_l \setminus \bd{\lambda}_l|$ is none other than the length of the chain $\mathcal{C}$. By a well-known
formula for crystal graphs relating weights and positions in the $i$-chains, since $\bd{\lambda}_l$ is the head of $\mathcal{C}$,
we have $k_i = ( \wt(\bd{\lambda}_l),\, \alpha_i ) = k$. Thus $\mathcal{C}$ is a $i$-chain of length $k_i$ with head $\bd{\lambda}_l$ and tail 
$\bd{\mu}_l$, which proves (iii). The divided powers
$e_i^{(k)}$, $f_i^{(k)} \in \Uq$ ($k \in \mN^*$) act on $\bd{\nu}_l \in \Fq$ as follows, with notation of Section 
\ref{paragraphe_3_actions_on_Lambdas}:
\begin{equation} \tag{$*$}
e_i^{(k)}.\bd{\nu}_l =
\ds \sum_{\bd{\kappa}_l \stackrel{i:k}{\longrightarrow} \bd{\nu}_l}
{q^{-M_i^{<}(\bd{\kappa}_l;\bd{\nu}_l;\bd{s}_l;n)} \,\bd{\kappa}_l} \qquad \mbox{and} \qquad
f_i^{(k)}.\bd{\nu}_l =
\ds \sum_{\bd{\nu}_l \stackrel{i:k}{\longrightarrow} \bd{\xi}_l} 
{q^{M_i^{>}(\bd{\nu}_l;\bd{\xi}_l;\bd{s}_l;n)} \,\bd{\xi}_l}.
\end{equation}
(For $k=1$, this is a part of Theorem \ref{thm_action_Uq_sur_Lambda_s}. The general case follows by induction on $k$; see \eg \cite{J} for a 
detailed proof.) Since $\bd{\mu}_l$ has no addable $i$-node and $\bd{\lambda}_l$ has no removable $i$-node, $\bd{\mu}_l \setminus \bd{\lambda}_l$ is 
the set of the addable $i$-nodes of $\bd{\lambda}_l$ and also the set of the removable $i$-nodes of $\bd{\mu}_l$, and this set has exactly $k_i$ 
elements. This proves (i), and this together with ($*$) proves (ii).
\end{proof} 

\begin{example} Take $n\!=\!3$, $l\!=\!2$, $\bd{s}_l\!=\!(1,2)$, $i\!=\!0$, $\bd{\lambda}_l=((2,2,1),(3,2))$ and 
$w=\wt(|\bd{\lambda}_l,\bd{s}_l \rangle)$. Then we have $\sigma_i(\bd{\lambda}_l)=((3,2,2),(3,3,1))$. \fini
\end{example}

\begin{remark} The proof of Proposition \ref{bij_scopes} shows that (with the assumptions and notation of this proposition) 
$\bd{\lambda}_l=\sigma_i^{-1}(\bd{\mu}_l)$ is the multi-partition obtained by removing to $\bd{\mu}_l$ all its removable $i$-nodes. \fini
\end{remark}

\subsection{A first theorem of comparison}

Define a symmetric bilinear non-degenerate form $(.\,,.)$ on $\Fq$ by
\begin{equation}
(\bd{\lambda}_l,\bd{\mu}_l)=q^{\| \bd{\lambda}_l \|} \delta_{\bd{\lambda}_l,\bd{\mu}_l} \qquad (\bd{\lambda}_l,\,\bd{\mu}_l \in \Pi^l), 
\end{equation}    
where we put
\begin{equation}
\| \bd{\lambda}_l \|:=(\wt(\bd{\lambda}_l),\wt(\bd{\lambda}_l))/2 \qquad (\bd{\lambda}_l \in \Pi^l).
\end{equation}
This form enjoys the following property:

\begin{lemma} \label{lemma_ps_Fqsl}
For $u$, $v \in \Fq$, $0 \leq i \leq n-1$, we have 
$$(t_i.u \,, v)=(u\, , t_i.v) \qquad \mbox{and} \qquad (e_i.u\,,v)=(u\,,f_i.v).$$
\end{lemma}

\begin{proof} Identical to the proof of \cite[Prop. 8.1]{LLT}. 
\end{proof}

Let $\set{G^{*}(\bd{\lambda}_l) \mid \bd{\lambda}_l \in \Pi^l}$ denote the adjoint basis of 
$\set{G^{+}(\bd{\lambda}_l) \mid \bd{\lambda}_l \in \Pi^l}$ with respect to the form $(.\,,.)$. Since the basis 
$\set{G^{+}(\bd{\lambda}_l) \mid \bd{\lambda}_l \in \Pi^l}$ is a lower global crystal basis of $\Fq$ in the sense of \cite{Kas1}, it follows by Lemma 
\ref{lemma_ps_Fqsl} and \cite[Prop. 3.2.2]{Kas1} that $\set{G^{*}(\bd{\lambda}_l) \mid \bd{\lambda}_l \in \Pi^l}$ is an upper global crystal basis of 
$\Fq$. \\

We are now ready to prove the following result, which is a generalization to higher-level Fock spaces of \cite[Thm. 20]{LM}.

\begin{thm} \label{thm_comparaison_bc_1} Let $\bd{s}_l \in \mZ^l(s)$, $w \in \mathcal{P}(\Fq)$ and $i \in \interv{0}{n-1}$ be such that 
$w + \alpha_i$ is not a weight of $\Fq$. Let $\sigma_i : \Pi^l(\bd{s}_l;w) \rightarrow \Pi^l(\bd{s}_l;\sigma_i.w)$ be the bijection defined by 
(\ref{def_sigma_i_2}). Then we have, for $\bd{\lambda}_l,\,\bd{\mu}_l \in \Pi^l(\bd{s}_l;w)$,
$$\mbox{\rm (i)} \quad \Delta^{+}_{\sigma_i(\bd{\lambda}_l),\sigma_i(\bd{\mu}_l) ;\, \bd{s}_l}(q) = 
\Delta^{+}_{\bd{\lambda}_l,\bd{\mu}_l ;\,\bd{s}_l}(q) \quad \mbox{and} \quad
\mbox{\rm (ii)} \quad 
\Delta^{-}_{\sigma_i(\bd{\lambda}_l),\sigma_i(\bd{\mu}_l);\, \bd{s}_l}(q)=\Delta^{-}_{\bd{\lambda}_l,\bd{\mu}_l;\, \bd{s}_l}(q).$$
As a consequence, the canonical bases of $\Fq \langle w \rangle$ and $\Fq \langle \sigma_i.w \rangle$ are similar in the sense of Definition
\ref{def_bc_semblables}.
\end{thm}

\begin{proof} Let us prove (i). Let $\bd{\mu}_l \in \Pi^l(\bd{s}_l;w)$.
Taking adjoint bases in (\ref{def_Delta_q}) yields  
$$q^{-\| \sigma_i(\bd{\mu}_l) \|} \sigma_i(\bd{\mu}_l) = \ds \sum_{\bd{\nu}_l \in \Pi^l(\bd{s}_l;\sigma_i.w)}
{\Delta^{+}_{\sigma_i(\bd{\mu}_l),\bd{\nu}_l ;\, \bd{s}_l}(q) \; G^{*}(\bd{\nu}_l)}.$$
Since $\sigma_i : \Pi^l(\bd{s}_l;w) \rightarrow \Pi^l(\bd{s}_l;\sigma_i.w)$ is a bijection, we can make in the sum above the reindexing
$\bd{\nu}_l = \sigma_i(\bd{\lambda}_l)$. If we now apply $e_i^{(k_i)}$ with $k_i:=(w,\alpha_i)$ to both hand-sides of this equality, we get  
\begin{equation}
\tag{$*$}
q^{-\| \sigma_i(\bd{\mu}_l) \|} e_i^{(k_i)}.\sigma_i(\bd{\mu}_l) = \ds \sum_{\bd{\lambda}_l \in \Pi^l(\bd{s}_l;w)}
{\Delta^{+}_{\sigma_i(\bd{\mu}_l),\sigma_i(\bd{\lambda}_l) ;\, \bd{s}_l}(q) \; e_i^{(k_i)}.G^{*}(\sigma_i(\bd{\lambda}_l))}.
\end{equation}
Note that
\begin{equation}
\tag{$**$}
\| \sigma_i(\bd{\mu}_l) \| = \big(\sigma_i.\wt(\bd{\mu}_l) \,,\, \sigma_i.\wt(\bd{\mu}_l) \big)/2 =
\big(\wt(\bd{\mu}_l) \,,\, \wt(\bd{\mu}_l) \big)/2 = \| \bd{\mu}_l \|.
\end{equation}
\noindent By Proposition \ref{bij_scopes} (ii), we have $e_i^{(k_i)}.\sigma_i(\bd{\mu}_l) = \bd{\mu}_l$. 
Now let $\bd{\lambda}_l \in \Pi^l(\bd{s}_l;w)$. Since $w + \alpha_i$ is not a weight of $\Fq$, we have
$e_i.\bd{\lambda}_l=0$, whence $\eiup.\bd{\lambda}_l=0$. Moreover, again by Proposition \ref{bij_scopes}, we have  
$\sigma_i(\bd{\lambda}_l)=f_i^{(k_i)}.\bd{\lambda}_l$. Therefore, by Lemma \ref{lemma_Kashiwara_operators}, we have
$$\bd{\lambda}_l = (\eilow)^{k_i}.(f_i^{(k_i)}.\bd{\lambda}_l) = (\eiup)^{k_i}.(f_i^{(k_i)}.\bd{\lambda}_l) =
(\eiup)^{k_i}.(\sigma_i(\bd{\lambda}_l)),$$
whence $(\eiup)^{k_i+1}.\sigma_i(\bd{\lambda}_l)=0$. Since $\set{G^{*}(\bd{\lambda}_l) \mid \bd{\lambda}_l \in \Pi^l}$ is an upper global crystal 
basis of $\Fq$, \cite[Lemma 5.1.1 (ii)]{Kas1} then implies
\begin{equation}
\tag{$***$}
e_i^{(k_i)}.G^{*}(\sigma_i(\bd{\lambda}_l))=G^{*} \big((\eiup)^{k_i}.(\sigma_i(\bd{\lambda}_l)) \big)= G^{*}(\bd{\lambda}_l).
\end{equation}
Combining ($*$), ($**$) and ($***$) we get

$$q^{-\| \bd{\mu}_l \|} \bd{\mu}_l = \ds \sum_{\bd{\lambda}_l \in \Pi^l(\bd{s}_l;w)}
{\Delta^{+}_{\sigma_i(\bd{\mu}_l),\sigma_i(\bd{\lambda}_l) ;\, \bd{s}_l}(q) \; G^{*}(\bd{\lambda}_l)}.$$
Since this is valid for any $\bd{\mu}_l \in \Pi^l(\bd{s}_l;w)$, we get the claimed formula by taking again adjoint bases. \\

Let us now prove (ii). Let $w' \in \mathcal{P}(\Fq)$. If we know the basis 
$\set{G^{+}(\bd{\lambda}_l) \mid \bd{\lambda}_l \in \Pi^l(\bd{s}_l;w')}$, then we can compute the involution of 
$\Fq \langle w' \rangle$ by solving a unitriangular system. Since the canonical basis 
$\set{G^{-}(\bd{\lambda}_l) \mid \bd{\lambda}_l \in \Pi^l(\bd{s}_l;w')}$ is uniquely determined by the involution of $\Fq \langle w' \rangle$,
the basis $\set{G^{-}(\bd{\lambda}_l) \mid \bd{\lambda}_l \in \Pi^l(\bd{s}_l;w')}$ is uniquely determined by the basis 
$\set{G^{+}(\bd{\lambda}_l) \mid \bd{\lambda}_l \in \Pi^l(\bd{s}_l;w')}$ (and conversely). Thus (i) implies (ii). (For a different proof
of this fact, one can also apply \cite[Thm 5.15]{U}.)
\end{proof}

\begin{example} \label{ex_comparaison_bc_1}
Take $n=3$, $l=2$, $\bd{s}_l=(1,0)$, $w=\wt(|\bd{\emptyset}_l, \bd{s}_l \rangle)-(2\alpha_0+3\alpha_1+\alpha_2)$ and
$i=2$. One can easily check that $w + \alpha_i$ is not a weight of $\Fq$. The elements of $\Pi^l(\bd{s}_l;w)$ are 
$$\scriptsize \begin{array}{l} {\big((1), (5) \big),\, \big((4), (2) \big),\, \big((4, 2), \emptyset \big),\, \big((1), (2, 2, 1) \big),\, 
\big((2, 2), (2) \big),\, \big((1, 1), (2, 1, 1) \big),\, \big((1, 1, 1, 1), (2) \big),\, \big((1), (2, 1, 1, 1) \big)} \, ,  
\end{array}$$
and their respective images by the map $\sigma_i$ are
$$\hspace*{-5mm} 
\scriptsize \begin{array}{l} {\big((2), (6,1) \big),\, \big((5), (3,1) \big),\, \big((5,3,1), \emptyset \big),\, \big((2), (3, 2, 2) \big),\, 
\big((2, 2,1), (3,1) \big),\, \big((2,1,1), (3, 1, 1) \big),\, \big((2, 1, 1, 1), (3,1) \big),\, \big((2), (3, 1, 1, 1,1) \big)} \, .
\end{array}$$
With obvious notation, the transition matrices of the canonical bases of $\Fq \langle w \rangle$ are

$$ \Delta_{\bd{s}_l}^+ \langle w \rangle (q) = \begin{array}{ll}
 \left( \begin{array}{cccccccc}
1& .& .& .& .& .& .& . \\
q& 1& .& .& .& .& .& . \\ 
0& q& 1& .& .& .& .& . \\ 
q& 0& 0& 1& .& .& .& . \\ 
q& q^2& q& 0& 1& .& .& . \\ 
q^2& 0& 0& q& q& 1& .& . \\
0& 0& q& 0& q^2& q& 1& . \\
0& 0& 0& q^2& 0& q& 0& 1
\end{array} \right) 
& \hspace{-3mm}
 \begin{array}{l}
 \big((1), (5) \big) \\
 \big((4), (2) \big) \\ 
 \big((4, 2), \emptyset \big) \\
 \big((1), (2, 2, 1) \big) \\
 \big((2, 2), (2) \big) \\ 
 \big((1, 1), (2, 1, 1) \big) \\ 
 \big((1, 1, 1, 1), (2) \big) \\ 
 \big((1), (2, 1, 1, 1) \big)  
 \end{array} \qquad \mbox{and}
\end{array}$$

$$ \Delta_{\bd{s}_l}^- \langle w \rangle (q) = \begin{array}{ll}
 \left( \begin{array}{cccccccc}
1        &  . 	     &  .  &  .  &  .  &  .  &  .  & . \\
 q^{-1}  &   1 	  	 &   . &   . &  .  &  .  &  .  & . \\
q^{-2}   & - q^{-1}  &   1 &  .  &  .  &  .  &  .  & . \\
- q^{-1} &   0 		 &   0 &   1 &  .  &  .  &  .  & . \\
0 		 &   0 		 & - q^{-1} &  0 &  1 &  .  &  .  &  .  \\
q^{-2}& - q^{-1} & q^{-2} & - q^{-1} & - q^{-1} & 1 & . & . \\
- q^{-3} &  q^{-2} &  0 & q^{-2} & 0 & - q^{-1} & 1 & . \\
0 &  q^{-2} &  - q^{-3} &  0 & q^{-2} & - q^{-1} & 0 & 1
\end{array} \right) 
& \hspace{-3mm}
 \begin{array}{l}
 \big((1), (5) \big) \\
 \big((4), (2) \big) \\ 
 \big((4, 2), \emptyset \big) \\
 \big((1), (2, 2, 1) \big) \\
 \big((2, 2), (2) \big) \\ 
 \big((1, 1), (2, 1, 1) \big) \\ 
 \big((1, 1, 1, 1), (2) \big) \\ 
 \big((1), (2, 1, 1, 1) \big)  
 \end{array} \, .
\end{array}$$

\bigskip

\noindent In the same way, the transition matrices of the canonical bases of $\Fq \langle \sigma_i.w \rangle$ are
$$ \Delta_{\bd{s}_l}^+ \langle \sigma_i.w \rangle (q) = \begin{array}{ll}
 \left( \begin{array}{cccccccc}
1& .& .& .& .& .& .& . \\
q& 1& .& .& .& .& .& . \\ 
0& q& 1& .& .& .& .& . \\ 
q& 0& 0& 1& .& .& .& . \\ 
q& q^2& q& 0& 1& .& .& . \\ 
q^2& 0& 0& q& q& 1& .& . \\
0& 0& q& 0& q^2& q& 1& . \\
0& 0& 0& q^2& 0& q& 0& 1
\end{array} \right) 
& \hspace{-3mm}
 \begin{array}{l}
\big((2), (6,1) \big) \\ 
\big((5), (3,1) \big) \\ 
\big((5,3,1), \emptyset \big) \\ 
\big((2), (3, 2, 2) \big) \\ 
\big((2, 2,1), (3,1) \big) \\ 
\big((2,1,1), (3, 1, 1) \big) \\ 
\big((2, 1, 1, 1), (3,1) \big) \\ 
\big((2), (3, 1, 1, 1,1) \big)  
 \end{array} \qquad \mbox{and}
\end{array}$$

$$ \Delta_{\bd{s}_l}^- \langle \sigma_i.w \rangle (q) = \begin{array}{ll}
 \left( \begin{array}{cccccccc}
1        &  .  &  .  &  .  &  .  &   . &   . &  . \\
 q^{-1}  &   1 &   . &   . &  .  &   . &   . &  . \\
q^{-2}   & - q^{-1} &   1 &   . &   . &  . &  . &  . \\
- q^{-1} &   0 &   0 &   1 &   . &   . &   . &   . \\
0 &   0 &   - q^{-1} &   0 &   1 &   . &   . &   . \\
q^{-2}& - q^{-1} & q^{-2} & - q^{-1} & - q^{-1} & 1 & . & . \\
- q^{-3} &  q^{-2} &  0 & q^{-2} & 0 & - q^{-1} & 1 & . \\
0 &  q^{-2} &  - q^{-3} &  0 & q^{-2} & - q^{-1} & 0 & 1
\end{array} \right) 
& \hspace{-3mm}
 \begin{array}{l}
\big((2), (6,1) \big) \\ 
\big((5), (3,1) \big) \\ 
\big((5,3,1), \emptyset \big) \\ 
\big((2), (3, 2, 2) \big) \\ 
\big((2, 2,1), (3,1) \big) \\ 
\big((2,1,1), (3, 1, 1) \big) \\ 
\big((2, 1, 1, 1), (3,1) \big) \\ 
\big((2), (3, 1, 1, 1,1) \big)  
 \end{array} \, ,
\end{array}$$

\noindent in agreement with Theorem \ref{thm_comparaison_bc_1}. \fini
\end{example}

\section{Comparison of canonical bases of $\Fq \langle w \rangle$ for $w$ in a given coset in $\mathcal{P}(\Lambda^s)/\mZ \delta$
and different multi-charges $\bd{s}_l$}

\begin{notation} \label{notation_thm_comparaison_bc_2}
For $\bd{a}_l$, $\bd{b}_l \in \mZ^l(s)$, introduce the shorter notation   
\begin{equation} \label{notation_d_sl_tl}
d(\bd{a}_l,\bd{b}_l):=\Delta(\bd{a}_l,n)- \Delta(\bd{b}_l,n) \in \mZ.
\end{equation}
Throughout this section we fix $\bd{s}_l \in \mZ^l(s)$, $w \in \mathcal{P}(\Fq)$ and $i \in \interv{0}{l-1}$. Finally, let
$(\bd{s}_n,\dot{w}) \in \mZ^n(s) \times \dot{\mathcal{P}}(\Lambda^s)$ be the pair such that
$\Fq \langle w \rangle = \Fp \langle \dot{w} \rangle$ (see Proposition \ref{prop_bij_pds_charges}). \fini
\end{notation}

\subsection{A second theorem of comparison}
Recall that $\dot{W}_l = \langle \dot{\sigma}_0,\ldots,\dot{\sigma}_{l-1} \rangle \cong \widetilde{\mathfrak{S}}_l$ is the Weyl group of $\Up$.

\begin{definition} \label{def_dot_sigma_i}
Keep Notation \ref{notation_thm_comparaison_bc_2}. By analogy with Definition \ref{def_sigma_i}, define a bijection 
\begin{equation} \label{eq_def_w_sl}
\dot{\sigma}_i : \Pi^n(\bd{s}_n;\dot{w}) \rightarrow \Pi^n(\bd{s}_n;\dot{\sigma}_i.\dot{w})
\end{equation}
which enjoys similar properties as the bijections $\sigma_j$ from (\ref{def_sigma_i_2}). Since 
$\Fq \langle w \rangle = \Fp \langle \dot{w} \rangle$, we have a bijection
between the standard basis of $\Fq \langle w \rangle$ (as a subspace of $\Fq$) and the standard basis of $\Fp \langle \dot{w} \rangle$
(as a subspace of $\Fp$). We thus have a bijection 
$\Pi^l(\bd{s}_l;w) \stackrel{\sim}{\longrightarrow} \Pi^n(\bd{s}_n;\dot{w})$. Put $\bd{t}_l:=\dot{\sigma}_i.\bd{s}_l$. 
The same argument gives, by Lemma \ref{lemma_action_gpe_Weyl_on_Fock_spaces}, a bijection  
$\Pi^l \big(\bd{t}_l;w+d(\bd{s}_l,\bd{t}_l)\,\delta \big) \stackrel{\sim}{\longrightarrow} \Pi^n(\bd{s}_n;\dot{\sigma}_i.\dot{w})$. We therefore have the 
following commutative diagram of bijections, in which the dashed arrow will still be denoted by $\dot{\sigma}_i$. 
$$\begin{array}{ccc}
\Pi^n(\bd{s}_n;\dot{w}) & \stackrel{\dot{\sigma}_i}{\longrightarrow} & \Pi^n(\bd{s}_n;\dot{\sigma}_i.\dot{w}) \\
\uparrow & & \uparrow \\
\Pi^l(\bd{s}_l;w) & \stackrel{\sim}{\dasharrow} & \Pi^l \big(\dot{\sigma}_i.\bd{s}_l\,;\,w+d(\bd{s}_l,\dot{\sigma}_i.\bd{s}_l) \, \delta \big).
\end{array}$$ \fini
\end{definition}

\begin{example} \label{ex_bij_sigma_i_pt}
Take $n=2$, $l=3$, $\bd{s}_l=(0,2,-1)$ and $i=2$. Note that $d(\bd{s}_l,\dot{\sigma}_i.\bd{s}_l)=0$ in this case. Take
$\bd{\lambda}_l=\big(\emptyset,(2),\emptyset \big) \in \Pi^l(\bd{s}_l;w)$, where 
$w:=\wt(|\bd{\emptyset}_l, \bd{s}_l \rangle)-(\alpha_0+\alpha_1)$. By Proposition \ref{prop_bij_pds_charges}, we have 
$\dot{w}=\dot{\wt}(|\bd{\emptyset}_n, \bd{s}_n \rangle^{\bullet})-(2 \dot{\alpha}_0 + 3\dot{\alpha}_1 + \dot{\alpha}_2)$
with $\bd{s}_n:=(1,0)$. We have $|\bd{\lambda}_l, \bd{s}_l \rangle = |\bd{\lambda}_n, \bd{s}_n \rangle^{\bullet}$ with
$\bd{\lambda}_n:=\big( (1),(5) \big)$. One computes $\bd{\mu}_n:=\dot{\sigma}_i(\bd{\lambda}_n)=\big( (2),(6,1) \big)$ (see Example 
\ref{ex_comparaison_bc_1}). Let $\bd{t}_l:=\dot{\sigma}_i.\bd{s}_l=(0,-1,2)$. Then $\dot{\sigma}_i(\bd{\lambda}_l)$ is the $l$-multi-partition such 
that $|\dot{\sigma}_i(\bd{\lambda}_l), \bd{t}_l \rangle = |\bd{\mu}_n, \bd{s}_n \rangle^{\bullet}$, namely 
$\dot{\sigma}_i(\bd{\lambda}_l)= \big(\emptyset,\emptyset,(2) \big)$. \fini
\end{example}

Since $\set{G^{-}(\bd{\lambda}_n, \bd{s}_n)^{\bullet} \mid \bd{\lambda}_n \in \Pi^n}$ is a lower global crystal basis of $\Fp$,
one can prove for the Fock space $\Fp$ an analogue of Theorem \ref{thm_comparaison_bc_1}. Rephrasing this result in terms of the indexation
$\bd{\lambda}_l$ leads to the following result.

\begin{thm} \label{thm_comparaison_bc_2}
Keep Notation \ref{notation_thm_comparaison_bc_2} and assume that $\dot{w}+\dot{\alpha}_i$ is not a weight of $\Fp$. Let 
$\dot{\sigma}_i : \Pi^l(\bd{s}_l;w) \rightarrow \Pi^l \big(\dot{\sigma}_i.\bd{s}_l\,;\,w+d(\bd{s}_l,\dot{\sigma}_i.\bd{s}_l) \, \delta \big)$ be the 
bijection from Definition \ref{def_dot_sigma_i}. Then we have, for 
$\bd{\lambda}_l,\,\bd{\mu}_l \in \Pi^l(\bd{s}_l;w)$, $\epsilon=\pm 1$:
$$\Delta^{\epsilon}_{\dot{\sigma}_i(\bd{\lambda}_l),\dot{\sigma}_i(\bd{\mu}_l) ;\, \dot{\sigma}_i.\bd{s}_l}(q) = 
\Delta^{\epsilon}_{\bd{\lambda}_l,\bd{\mu}_l ;\, \bd{s}_l}(q).$$ 
\end{thm}

\begin{proof} Apply the analogue for $\Fp$ of Theorem \ref{thm_comparaison_bc_1} mentioned above, then use Lemma 
\ref{lemma_action_gpe_Weyl_on_Fock_spaces} and (\ref{equation_lien_Delta_Deltapoint}).
\end{proof}

\begin{example} \label{ex_comparaison_bc_2} \emph{(see Exemples \ref{ex_comparaison_bc_1} and \ref{ex_bij_sigma_i_pt})}
Take $n$, $l$, $\bd{s}_l$, $w$ and $i$ as in Example \ref{ex_bij_sigma_i_pt} (namely,
$n:=2$, $l:=3$, $\bd{s}_l:=(0,2,-1)$, $w:=\wt(|\bd{\emptyset}_l, \bd{s}_l \rangle)-(\alpha_0+\alpha_1)$ and $i:=2$).
Note that $\dot{w}+\dot{\alpha}_i \notin \dot{\mathcal{P}}(\Fp)$ (see Example \ref{ex_comparaison_bc_1}). The elements of $\Pi^l(\bd{s}_l;w)$ are 
$$\scriptsize \begin{array}{l} {\big( \emptyset, (2), \emptyset \big) ,\,\big( \emptyset, \emptyset, (1,1) \big) ,\,
\big( \emptyset, (1) , (1) \big) ,\,\big( \emptyset, \emptyset, (2) \big) ,\,\big( (2), \emptyset, \emptyset \big) ,\,
\big( (1), \emptyset, (1) \big) ,\,\big( (1,1), \emptyset, \emptyset \big) ,\,\big( \emptyset, \emptyset, (1,1) \big) },
\end{array}$$
and their respective images by the map $\dot{\sigma_i}$ are

$$\scriptsize \begin{array}{l} {\big( \emptyset, \emptyset, (2) \big) ,\,\big( \emptyset, \emptyset, (1,1) \big) ,\,
\big( \emptyset, (1) , (1) \big) ,\,\big( (2), \emptyset, \emptyset \big) ,\,\big( \emptyset,(2), \emptyset \big) ,\,
\big( (1), (1), \emptyset \big) ,\,\big( (1,1), \emptyset, \emptyset \big) ,\,\big( \emptyset, (1,1), \emptyset \big) },
\end{array}$$

\bigskip

\noindent On the one hand, the transition matrices of the canonical bases of $\Fq \langle w \rangle$ are

$$ \Delta_{\bd{s}_l}^+ \langle w \rangle (q) = \begin{array}{ll}
\left( \begin{array}{cccccccc}
1& .& .& .& .& .& .& . \\ 
q& 1& .& .& .& .& .& . \\ 
q^2& q& 1& .& .& .& .& . \\ 
0& 0& q& 1& .& .& .& . \\ 
q& 0& 0& 0& 1& .& .& . \\ 
q^2& q& q^2& q& q& 1& .& . \\ 
q^3& q^2& 0& 0& q^2& q& 1& . \\ 
0& q^2& q^3& q^2& 0& q& 0& 1 
\end{array} \right) 
& \hspace{-3mm}
\begin{array}{l}
\big( \emptyset, (2), \emptyset \big)  \\
\big( \emptyset, \emptyset, (1,1) \big)  \\
\big( \emptyset, (1) , (1) \big)  \\
\big( \emptyset, \emptyset, (2) \big)  \\
\big( (2), \emptyset, \emptyset \big)  \\
\big( (1), \emptyset, (1) \big)  \\
\big( (1,1), \emptyset, \emptyset \big)  \\
\big( \emptyset, \emptyset, (1,1) \big)
\end{array} \qquad \mbox{and}
\end{array}$$

$$ \Delta_{\bd{s}_l}^- \langle w \rangle (q) = \begin{array}{ll}
 \left( \begin{array}{cccccccc}
 
1& .& .& .& .& .& .& . \\ 
-q^{-1}& 1& .& .& .& .& .& . \\ 
0& -q^{-1}& 1& .& .& .& .& . \\ 
-q^{-1}& q^{-2}& -q^{-1}& 1& .& .& .& . \\ 
-q^{-1}& 0& 0& 0& 1& .& .& . \\ 
q^{-2}& 0& 0& -q^{-1}& -q^{-1}& 1& .& . \\ 
0& 0& -q^{-1}& q^{-2}& 0& -q^{-1}& 1& . \\ 
0& 0& 0& 0& q^{-2}& -q^{-1}& 0& 1
\end{array} \right) 
& \hspace{-3mm}
 \begin{array}{l}
\big( \emptyset, (2), \emptyset \big)  \\
\big( \emptyset, \emptyset, (1,1) \big)  \\
\big( \emptyset, (1) , (1) \big)  \\
\big( \emptyset, \emptyset, (2) \big)  \\
\big( (2), \emptyset, \emptyset \big)  \\
\big( (1), \emptyset, (1) \big)  \\
\big( (1,1), \emptyset, \emptyset \big)  \\
\big( \emptyset, \emptyset, (1,1) \big)  
\end{array} \, .
\end{array}$$
\noindent On the other hand, the transition matrices of the canonical bases of $\Fq \langle w \rangle$ are
$$ \Delta_{\dot{\sigma_i}.\bd{s}_l}^+ \langle w \rangle (q) = \begin{array}{ll}
 \left( \begin{array}{cccccccc}
1& .& .& .& .& .& .& . \\ 
q& 1& .& .& .& .& .& . \\ 
q^2& q& 1& .& .& .& .& . \\ 
0& 0& q& 1& .& .& .& . \\ 
q& 0& 0& 0& 1& .& .& . \\ 
q^2& q& q^2& q& q& 1& .& . \\ 
q^3& q^2& 0& 0& q^2& q& 1& . \\ 
0& q^2& q^3& q^2& 0& q& 0& 1 
\end{array} \right) 
& \hspace{-3mm}
 \begin{array}{l}
 \big( \emptyset, \emptyset, (2) \big)  \\
 \big( \emptyset, \emptyset, (1,1) \big)  \\
 \big( \emptyset, (1) , (1) \big)  \\
 \big( (2), \emptyset, \emptyset \big)  \\
 \big( \emptyset,(2), \emptyset \big)  \\
 \big( (1), (1), \emptyset \big)  \\
 \big( (1,1), \emptyset, \emptyset \big)  \\
 \big( \emptyset, (1,1), \emptyset \big) 
 \end{array} \qquad \mbox{and}
\end{array}$$

$$ \Delta_{\dot{\sigma_i}.\bd{s}_l}^- \langle w \rangle (q) = \begin{array}{ll}
 \left( \begin{array}{cccccccc}
1& .& .& .& .& .& .& . \\ 
-q^{-1}& 1& .& .& .& .& .& . \\ 
0& -q^{-1}& 1& .& .& .& .& . \\ 
-q^{-1}& q^{-2}& -q^{-1}& 1& .& .& .& . \\ 
-q^{-1}& 0& 0& 0& 1& .& .& . \\ 
q^{-2}& 0& 0& -q^{-1}& -q^{-1}& 1& .& . \\ 
0& 0& -q^{-1}& q^{-2}& 0& -q^{-1}& 1& . \\ 
0& 0& 0& 0& q^{-2}& -q^{-1}& 0& 1
\end{array} \right) 
& \hspace{-3mm}
 \begin{array}{l}
 \big( \emptyset, \emptyset, (2) \big)  \\
 \big( \emptyset, \emptyset, (1,1) \big)  \\
 \big( \emptyset, (1) , (1) \big)  \\
 \big( (2), \emptyset, \emptyset \big)  \\
 \big( \emptyset,(2), \emptyset \big)  \\
 \big( (1), (1), \emptyset \big)  \\
 \big( (1,1), \emptyset, \emptyset \big)  \\
 \big( \emptyset, (1,1), \emptyset \big) 
 \end{array} \, ,
\end{array}$$
\noindent in agreement with Theorem \ref{thm_comparaison_bc_2} (note that $d(\bd{s}_l,\dot{\sigma_i}.\bd{s}_l)=0$ in this case). \fini
\end{example}

\subsection{A sufficient condition for Theorem \ref{thm_comparaison_bc_2}} \label{paragraphe_CS_pour_thm_comparaison_bc_2}
Keep again Notation \ref{notation_thm_comparaison_bc_2}. We can apply Theorem \ref{thm_comparaison_bc_2} if 
$\dot{w}+\dot{\alpha}_i$ is not a weight of $\Fp$. By \cite[Prop. 3.6~(iv)]{Kac}, this condition holds if and only if 
$\dot{e}_i.(|\bd{\lambda}_l,\bd{s}_l \rangle)=0$ for all $\bd{\lambda}_l \in \Pi^l(\bd{s}_l;w)$.
We therefore have to check, for all $\bd{\lambda}_l \in \Pi^l(\bd{s}_l;w)$, whether
$\bd{\lambda}_n$ has a removable $i$-node, where $\bd{\lambda}_n \in \Pi^n$ is related to $\bd{\lambda}_l$ by 
$|\bd{\lambda}_l,\bd{s}_l \rangle = |\bd{\lambda}_n,\bd{s}_n \rangle^{\bullet}$. It is not very convenient to make such tests in
practice when the cardinality of $\Pi^l(\bd{s}_l;w)$ is large. We shall therefore give a sufficient condition on $\bd{s}_l$ and $w$ that ensures, without further computation,
that $\dot{w}+\dot{\alpha}_i$ is not a weight of $\Fp$. 

\begin{notation} \label{notation_Ni_w}
Let $\bd{s}_l \in \mZ^l(s)$, $w \in \mathcal{P}(\Fq)$ and $i \in \interv{0}{l-1}$. By (\ref{poidsl}),
the integer $N_i(\bd{\lambda}_l;\bd{s}_l;n)$ only depends on $\bd{s}_l$ and $w$, but not on $\bd{\lambda}_l \in \Pi^l(\bd{s}_l;w)$.
In the sequel this number will be denoted by $N_i(w;\bd{s}_l)$. \fini
\end{notation}

\begin{lemma} \label{lemma_Ni_w}  
Keep Notation \ref{notation_thm_comparaison_bc_2}. Then for all $\dot{\sigma} \in \dot{W}_l$, $0 \leq i \leq n-1$,
$w \in \mathcal{P}(\Fq)$, we
have $$N_i\big( w+d(\bd{s}_l,\dot{\sigma}.\bd{s}_l)\, \delta\,;\,\dot{\sigma}.\bd{s}_l\big)=N_i(w;\bd{s}_l).$$
\end{lemma}

\begin{proof} By (\ref{poidsl}) and the definition of the integers 
$N_j\big( w+d(\bd{s}_l,\dot{\sigma}.\bd{s}_l)\, \delta\,;\,\dot{\sigma}.\bd{s}_l\big)$ and $d(\bd{s}_l,\dot{\sigma}.\bd{s}_l)$, we have
$$\begin{array}{rcl}
w + d(\bd{s}_l,\dot{\sigma}.\bd{s}_l)\, \delta &=& \wt(|\bd{\emptyset}_l, \dot{\sigma}.\bd{s}_l \rangle) 
-\ds \sum_{j=0}^{n-1} N_j\big( w+d(\bd{s}_l,\dot{\sigma}.\bd{s}_l)\, \delta\,;\,\dot{\sigma}.\bd{s}_l\big) \, \alpha_j \\
 &=& \wt(|\bd{\emptyset}_l, \bd{s}_l \rangle) + d(\bd{s}_l,\dot{\sigma}.\bd{s}_l)\,\delta
-\ds \sum_{j=0}^{n-1} N_j\big( w+d(\bd{s}_l,\dot{\sigma}.\bd{s}_l)\, \delta\,;\,\dot{\sigma}.\bd{s}_l\big) \, \alpha_j,
\end{array}$$
whence $w = \wt(|\bd{\emptyset}_l, \bd{s}_l \rangle) 
-\ds \sum_{j=0}^{n-1} N_j\big( w+d(\bd{s}_l,\dot{\sigma}.\bd{s}_l)\, \delta\,;\,\dot{\sigma}.\bd{s}_l\big) \, \alpha_j.$
The lemma follows.
\end{proof}

\begin{lemma} \label{lemme_CS_pour_thm_comparaison_bc_2} 
Keep Notation \ref{notation_thm_comparaison_bc_2}. Assume that
$$s_i-s_{i+1} \geq n\big(N_0(w\,;\bd{s}_l)+1\big)$$
(where we put $s_0:=n+s_l$ if $i=0$). Then $\dot{w}+\dot{\alpha}_i$ is not a weight of $\Fp$. 
\end{lemma}

\begin{proof} We only sketch this proof and leave the details to the reader. Assume on the contrary that  
$\dot{w}+\dot{\alpha}_i$ is a weight of $\Fp$. Let $(\bd{\lambda}_l,\bd{t}_l) \in \Pi^l \times \mZ^l(s)$ be such that 
$\wtpt(|\bd{\lambda}_l,\bd{t}_l \rangle)=\dot{w}+\dot{\alpha}_i$. By (\ref{poidsptl}) we get
$$N_0(\bd{\lambda}_l;\bd{t}_l;n) = \Delta(\bd{s}_l,n) - \Delta(\bd{t}_l,n) + N_0(w\,;\bd{s}_l) - \delta_{i,0}$$
and $\bd{t}_l=(t_1,\ldots,t_l)$, where $t_0,\ldots,t_l$ are defined by
$$t_0= t_l+n, \quad t_i := s_i+1, \quad t_{i+1} := s_{i+1}-1  \quad \mbox{and} \quad t_j:=s_j \mbox{ if } j \notin \set{i,i+1,i+l}.$$ 
Moreover, we have by assumption $$N_0(w\,;\bd{s}_l)-\frac{1}{n}(s_i-s_{i+1}) \leq -1.$$
Combining these facts gives, after a little computation: $N_0(\bd{\lambda}_l;\bd{t}_l;n)<0$, which is absurd since
$N_0(\bd{\lambda}_l;\bd{t}_l;n)$ is the number of $0$-nodes of the multi-partition $\bd{\lambda}_l$.  
\end{proof}

\begin{remark} \label{rq_raffinement_lemme_CS_pour_thm_comparaison_bc_2}
The lower bound $s_i-s_{i+1} \geq n\big(N_0(w\,;\bd{s}_l)+1\big)$ from Lemma \ref{lemme_CS_pour_thm_comparaison_bc_2} is certainly not the 
best to ensure that $\dot{w}+\dot{\alpha}_i \notin \mathcal{P}(\Fp)$. We actually conjecture that the latter statement holds if 
$$s_i-s_{i+1} \geq N_0(w\,;\bd{s}_l)+\cdots+N_{n-1}(w\,;\bd{s}_l)$$ 
(this lower bound is in general better). \fini
\end{remark}

\subsection{A graph containing multi-charges conjugated under the action of $\dot{W}_l$} 

\begin{definition} \label{definition_Gamma_M}
Fix $\bd{r}_l \in A_{l,n}(s)$ and $M \in \mN^*$. (Recall that the set $A_{l,n}(s)$ defined by (\ref{definition_A_l_n_s}) is a fundamental domain 
for the action of $\dot{W}_l$.) For $\bd{s}_l=(s_1,\ldots,s_l) \in \dot{W}_l.\bd{r}_l$, $\bd{t}_l \in \dot{W}_l.\bd{r}_l$, $0 \leq i \leq l-1$, write 
$\bd{s}_l \stackrel{i}{\longrightarrow} \bd{t}_l$ if $\bd{t}_l = \dot{\sigma}_i.\bd{s}_l$ and $s_i-s_{i+1} \geq M$
(if $i=0$, we put $s_0:=n+s_l$). Let $\Gamma(M)$ be the graph containing $\dot{W}_l.\bd{r}_l$ as set of vertices and the arrows 
$\bd{s}_l \stackrel{i}{\longrightarrow} \bd{t}_l$ ($\bd{s}_l$, $\bd{t}_l \in \dot{W}_l.\bd{r}_l$, $0 \leq i \leq l-1$). \fini
\end{definition}

\begin{remark} \label{finitude_nb_cc_de_Gamma_M} 
Note that $\Gamma(1)$ is connected. More generally, we claim that $\Gamma(M)$ has finitely many connected components.  
For $\bd{s}_l \in \dot{W}_l.\bd{r}_l$, let $\dot{\sigma}(\bd{s}_l) \in \dot{W}_l$ be the element of minimal length such that 
$\bd{s}_l = \dot{\sigma}(\bd{s}_l).\bd{r}_l$. Now, for each connected component $C$ in $\Gamma(M)$, choose 
$\bd{s}_l(C) \in C$ such that $\ell(\bd{s}_l(C))$ is minimal in the set $\set{\ell(\bd{t}_l) \mid \bd{t}_l \in C}$. (We think that this determines 
$\bd{s}_l(C)$ in a unique way, but we do not have the proof for this fact.) One can then show easily that $\bd{s}_l(C)$ lies in the finite set 
$$\set{(t_1,\ldots,t_l) \in \mZ^l(s) \mid \forall \ 0 \leq i \leq l-1, \, t_{i+1}-t_i \leq M-1},$$ 
which proves the claim. \fini
\end{remark}

We now give a relation between the (connected components of) $\Gamma(M)$ and the comparison of canonical bases. We shall give an important
application of this result in Section \ref{paragraphe_comparaison_bc_cas_dominant}.			  

\begin{prop} \label{prop_comparaison_bc_2}  

Let $\bd{r}_l \in A_{l,n}(s)$ and $w \in \mathcal{P}(\Fq[\bd{r}_l])$. Put $M:=n(N_0(w\,;\bd{r}_l)+1)$. Let $\bd{s}_l$ and $\bd{t}_l$ be two multi-charges in the same 
connected component of $\,\Gamma(M)$ (in particular, $\bd{s}_l$ and $\bd{t}_l$ are in the $\dot{W}_l$-orbit of $\bd{r}_l$). Then 
the canonical bases of $\,\Fq \langle w+d(\bd{r}_l,\bd{s}_l)\, \delta \rangle$ and $\Fq[\bd{t}_l] \langle w+d(\bd{r}_l,\bd{t}_l)\, \delta \rangle$ are similar in the sense of 
Definition \ref{def_bc_semblables}.
\end{prop}   

\begin{proof} 
We may assume that $\bd{s}_l \stackrel{i}{\longrightarrow} \bd{t}_l$ with $i \in \interv{0}{l-1}$. With obvious notation, we have
by Lemma \ref{lemma_Ni_w} : $s_i-s_{i+1} \geq n(N_0(w\,;\bd{r}_l)+1)=n\big(N_0(w+d(\bd{r}_l,\bd{s}_l)\, \delta \,;\bd{s}_l)+1\big)$. 
We can therefore apply Lemma \ref{lemme_CS_pour_thm_comparaison_bc_2} and then Theorem \ref{thm_comparaison_bc_2} to conclude.
\end{proof}

With the notation above, Proposition \ref{prop_comparaison_bc_2} and Remark \ref{finitude_nb_cc_de_Gamma_M} show
that there are only finitely many similarity classes of canonical bases of $\Fq \langle w+d(\bd{r}_l,\bd{s}_l)\, \delta \rangle$,
where $\bd{s}_l$ ranges over the $\dot{W}_l$-orbit of $\bd{r}_l$ and $(\bd{r}_l,w)$ is fixed.

\section{Comparison of canonical bases for dominant multi-charges} \label{paragraphe_comparaison_bc_cas_dominant}
\begin{definition} \label{def_charge_dominante}
Let $M \in \mN$. We say that $(x_1,\ldots,x_N) \in \mR^N$ is \emph{$M$-dominant} if for all $1 \leq i \leq N-1$, we have
$$x_i-x_{i+1} \geq M.$$ \fini  
\end{definition}

Throughout Section \ref{paragraphe_comparaison_bc_cas_dominant}, we keep the following notation.

\subsection{Notation} \label{notation_paragraphe_comparaison_bc_cas_dominant}
\begin{itemize}

\item[*] Recall that $\mR^l(s)=\set{(x_1,\ldots,x_l) \in \mR^l \mid x_1+\cdots+x_l=s}$. The subset of $\mR^l(s)$ formed
by the $M$-dominant elements will be denoted by $\mathcal{C}_M$. \\

\item[*] The group $\dot{W}_l \cong \widetilde{\mathfrak{S}}_l$ is a semidirect product of the finite symmetric group $\mathfrak{S}_l$ and
an abelian group $\dot{\mathcal{Q}}$ which is free of rank $l-1$. More precisely, $\dot{\mathcal{Q}}$ is spanned by 
$\dot{\tau}_1,\ldots,\dot{\tau}_{l-1}$, where $\dot{\tau}_i$ $(1 \leq i \leq l-1)$ acts on 
$\mZ^l(s)$ by
\begin{equation}
\dot{\tau}_i.(s_1,\ldots,s_l)=(s_1,\ldots,s_{i-1},s_i+n,s_{i+1}-n,s_{i+2},\ldots,s_l) \quad \big((s_1,\ldots,s_l) \in \mZ^l(s) \big). 
\end{equation} 
(since $\dot{W}_l$ acts faithfully on $\mZ^l(s)$, this determines $\dot{\tau}_i$ completely). \\

\item[*] For $\bd{a}_l=(a_1,\ldots,a_l) \in \mZ^l$, put
\begin{equation}
\mathcal{L}_{\bd{a}_l}:= \set{ (s_1,\ldots,s_l) \in \mZ^l(s) \mid \forall \, 1 \leq i \leq l, \; s_i \equiv a_i \pmod n }.
\end{equation}
Note that $\dot{\mathcal{Q}}$ acts transitively on $\mathcal{L}_{\bd{a}_l}$; in particular, two elements
of $\mathcal{L}_{\bd{a}_l}$ lie in a same $\dot{W}_l$-orbit. \\ 

\item[*] For $\bd{s}_l$, $\bd{t}_l \in \mZ^l(s)$ and $M \in \mN$, write
\begin{equation}
\bd{s}_l \Zconnected{M} \bd{t}_l
\end{equation} 
if $\mathcal{L}_{\bd{s}_l}=\mathcal{L}_{\bd{t}_l}$, and
there exist $\bd{s}_l^{(0)},\ldots,\bd{s}_l^{(r)} \in \mathcal{L}_{\bd{s}_l} \cap \mathcal{C}_M$ such that 
$\bd{s}_l^{(0)}=\bd{s}_l$, $\bd{s}_l^{(r)}=\bd{t}_l$ and for all $1 \leq i \leq r$, we have 
$\bd{s}_l^{(i)} \in \set{\dot{\tau}_j^{\pm 1}.\bd{s}_l^{(i-1)} \mid 1 \leq j \leq l-1}$. In other words, put into a non-oriented graph all the 
elements of $\mZ^l(s)$ and draw an edge between two vertices if they are $\dot{W}_l$-conjugated to each other by a generator
of $\dot{\mathcal{Q}}$ or its inverse. Then $\bd{s}_l \Zconnected{M} \bd{t}_l$ if and only if there exists a path in this graph connecting 
$\bd{s}_l$ and $\bd{t}_l$ through $M$-dominant vertices (including $\bd{s}_l$ and $\bd{t}_l$). \\

\item[*] Let $\bd{r}_l=(r_1,\ldots,r_l) \in A_{l,n}(s)$ and $w \in \mathcal{P}(\Fq[\bd{r}_l])$. 
Recall the definition of the integers $d(\bd{r}_l,\bd{s}_l)$ $(\bd{s}_l \in \mZ^l(s))$ from 
(\ref{notation_d_sl_tl}) and $N_i(w;\bd{r}_l)$ $(0 \leq i \leq n-1)$ from Notation \ref{notation_Ni_w}.
\end{itemize}

\subsection{A third theorem of comparison}

The goal of Section \ref{paragraphe_comparaison_bc_cas_dominant} is to prove the following theorem. 

\begin{thm} \label{thm_comparaison_bc_cas_dominant}
 
Keep Notation \ref{notation_paragraphe_comparaison_bc_cas_dominant}. Then there exists $N \in \mN$ 
(which only depends on $n$, $l$ and $N_0(w;\bd{r}_l)$) such that for all $N$-dominant multi-charges 
$\bd{s}_l$, $\bd{t}_l \in \dot{W}_l.\bd{r}_l$ with $\mathcal{L}_{\bd{s}_l}=\mathcal{L}_{\bd{t}_l}$, the canonical bases 
of $\Fq \langle w+d(\bd{r}_l,\bd{s}_l)\,\delta \rangle$ and 
$\Fq[\bd{t}_l] \langle w+d(\bd{r}_l,\bd{t}_l)\,\delta \rangle$ are similar. \cqfd
\end{thm}

\begin{remark} \label{rq_comparaison_bc_cas_dominant}
The proof of Theorem \ref{notation_paragraphe_comparaison_bc_cas_dominant} will provide an integer $N$ of the form
$$N=n N_0(w;\bd{r}_l)+c,$$ where $c$ is a constant that can be explicitly calculated.
More precisely, by Remark \ref{rq_constante_prop2_thm_comparaison_bc_cas_dominant}
we can take $c \geq n(l^2+l+3)$. However, this value of $N$
is probably not optimal. Indeed, according to Remark \ref{rq_raffinement_lemme_CS_pour_thm_comparaison_bc_2} and 
explicit calculations of canonical bases, we conjecture that Theorem \ref{thm_comparaison_bc_cas_dominant} holds if we replace
this $N$ by $$N':=N_0(w;\bd{r}_l)+\cdots+N_{n-1}(w;\bd{r}_l)$$
(the latter lower bound is in general better). \fini
\end{remark}

\begin{example} \label{ex_thm_comparaison_bc_cas_dominant}
Take $n=3$, $l=2$, $\bd{r}_l=(1,0)$ and $w=\wt(|\bd{\emptyset}_l, \bd{r}_l \rangle) - (\alpha_0+\alpha_1+\alpha_2)$.
With notation from Theorem \ref{thm_comparaison_bc_cas_dominant} and Remark \ref{rq_comparaison_bc_cas_dominant}, we 
can take $N=30$ and $N'=3$. Note that all the multi-charges $\bd{s}_l^{(k)}:=(3k+1,-3k)$ $(k \in \mZ)$ are in the $\dot{W}_l$-orbit of $\bd{r}_l$
and they have the same pair of residues modulo $n$. Put $w_k:=w+d(\bd{r}_l,\bd{s}_l^{(k)})\,\delta$ ($k \in \mZ$).
By Theorem \ref{thm_comparaison_bc_cas_dominant}, the canonical bases of $\Fq[\bd{s}_l^{(k)}] \langle w_k \rangle$, $k \geq 5$
are pairwise similar. By Remark \ref{rq_comparaison_bc_cas_dominant}, the canonical bases of $\Fq[\bd{s}_l^{(k)}] \langle w_k \rangle$, $k \geq 1$
should be pairwise similar, which can be actually checked by an explicit calculus. Namely, the transition matrices $\Delta^{\epsilon}_k(q)$ 
of the canonical bases of $\Fq[\bd{s}_l^{(k)}] \langle w_k \rangle$ ($k \geq 1$, $\epsilon = \pm 1$) are

$$ \Delta_k^+ (q) = \begin{array}{ll}
\left( \begin{array}{cccccccc}
1 & . & . & . & . & . & . & . \\ 
q & 1 & . & . & . & . & . & . \\ 
q^2 & q & 1 & . & . & . & . & . \\ 
0 & q & 0 & 1 & . & . & . & . \\ 
0 & q^2 & q & q & 1 & . & . & . \\ 
0 & 0 & q & 0 & 0 & 1 & . & . \\ 
0 & 0 & q^2 & 0 & q & q & 1 & . \\ 
0 & 0 & 0 & q & q^2 & 0 & q & 1
\end{array} \right) 
& \hspace{-3mm}
\begin{array}{l}
\big( (3), \emptyset \big) \\ 
\big( (2, 1), \emptyset \big) \\ 
\big( (2), (1) \big) \\ 
\big( (1, 1, 1), \emptyset \big) \\ 
\big( (1), (1, 1) \big) \\ 
\big( \emptyset, (3) \big) \\ 
\big( \emptyset, (2, 1) \big) \\ 
\big( \emptyset, (1, 1, 1) \big)
\end{array} \qquad \mbox{and}
\end{array}$$

$$ \Delta_k^- (q) = \begin{array}{ll}
\left( \begin{array}{cccccccc}
1 & . & . & . & . & . & . & . \\ 
-q^{-1} & 1 & . & . & . & . & . & . \\ 
0 & -q^{-1} & 1 & . & . & . & . & . \\ 
q^{-2} & -q^{-1} & 0 & 1 & . & . & . & . \\ 
0 & q^{-2} & -q^{-1} & -q^{-1} & 1 & . & . & . \\ 
-q^{-1} & q^{-2} & -q^{-1} & 0 & 0 & 1 & . & . \\ 
q^{-2} & -q^{-3} & q^{-2} & q^{-2} & -q^{-1} & -q^{-1} & 1 & . \\ 
-q^{-3} & 0 & 0 & 0 & 0 & q^{-2} & -q^{-1} & 1
\end{array} \right) 
& \hspace{-3mm}
\begin{array}{l}
\big( (3), \emptyset \big) \\ 
\big( (2, 1), \emptyset \big) \\ 
\big( (2), (1) \big) \\ 
\big( (1, 1, 1), \emptyset \big) \\ 
\big( (1), (1, 1) \big) \\ 
\big( \emptyset, (3) \big) \\ 
\big( \emptyset, (2, 1) \big) \\ 
\big( \emptyset, (1, 1, 1) \big)
\end{array} 
\end{array}.$$

One can check moreover that

$$ \Delta_0^+ (q) = \begin{array}{ll}
\left( \begin{array}{cccccccc}
1 & . & . & . & . & . & . & . \\ 
0 & 1 & . & . & . & . & . & . \\ 
0 & q & 1 & . & . & . & . & . \\ 
q & q & 0 & 1 & . & . & . & . \\ 
0 & q^2 & q & q & 1 & . & . & . \\ 
0 & 0 & q^2 & 0 & q & 1 & . & . \\ 
q^2 & 0 & 0 & q & 0 & 0 & 1 & . \\ 
0 & 0 & 0 & q^2 & q & 0 & q & 1
\end{array} \right) 
& \hspace{-3mm}
\begin{array}{l}
\big( (3), \emptyset \big) \\ 
\big( (2, 1), \emptyset \big) \\ 
\big( (2), (1) \big) \\ 
\big( (1, 1, 1), \emptyset \big) \\ 
\big( (1), (1, 1) \big) \\ 
\big( \emptyset, (3) \big) \\ 
\big( \emptyset, (2, 1) \big) \\ 
\big( \emptyset, (1, 1, 1) \big)
\end{array}.
\end{array}$$

\bigskip

\noindent Note that $\Delta_1^+ (q)$ a $22$ nonzero entries, whereas $\Delta_0^+ (q)$ has only $21$ nonzero entries.
As a consequence, the canonical bases of $\Fq[\bd{s}_l^{(k)}] \langle w_k \rangle$ for $k = 0$ and $k = 1$ are not similar. \fini
\end{example}

\bigskip

The proof of Theorem \ref{thm_comparaison_bc_cas_dominant} relies on the two following propositions.

\begin{prop} \label{prop1_thm_comparaison_bc_cas_dominant} 
Keep the notation from Section \ref{notation_paragraphe_comparaison_bc_cas_dominant}. Let $\bd{s}_l \in \dot{W}_l.\bd{r}_l$ be
an $M$-dominant multi-charge with $M:=n(N_0(w;\bd{r}_l)+2)$. Let $1 \leq i \leq l-1$ and $\bd{t}_l := \dot{\tau}_i.\bd{s}_l$. Then the canonical bases of
$\Fq \langle w+d(\bd{r}_l,\bd{s}_l)\,\delta \rangle$ and $\Fq[\bd{t}_l] \langle w+d(\bd{r}_l,\bd{t}_l)\,\delta \rangle$ are similar.
\end{prop}

\begin{proof} We give the proof for $2 \leq i \leq l-3$ (the proof for $i=1$, $i=l-2$ and $i=l-1$ is similar).
One easily checks that
$$\dot{\tau}_i = \dot{\sigma}_{i-1} \dot{\sigma}_{i-2} \cdots \dot{\sigma}_1 \dot{\sigma}_0 \dot{\sigma}_{l-1} \cdots \dot{\sigma}_{i+2} 
\dot{\sigma}_{i+1} \dot{\sigma}_{i+2} \cdots \dot{\sigma}_{l-1} \dot{\sigma}_0 \dot{\sigma}_1 \cdots \dot{\sigma}_{i}.$$
Let $0 \leq k \leq 2l-2$. Denote by $\dot{\tau}_i[k]$ the right factor of length $k$ in this word (we thus have
$\dot{\tau}_i[0]=\mbox{id}$, $\dot{\tau}_i[1]=\dot{\sigma}_{i}$, $\dot{\tau}_i[2]=\dot{\sigma}_{i-1} \dot{\sigma}_{i}$ and so on). 
Put $$\bd{s}_l^{(k)}=(s^{(k)}_1,\ldots,s^{(k)}_l):= \dot{\tau}_i[k].\bd{s}_l.$$
For $1 \leq k \leq 2l-2$, let $i_k \in \interv{0}{l-1}$ be such that $\dot{\tau}_i[k]=\dot{\sigma}_{i_k} \dot{\tau}_i[k-1]$. Let
now $0 \leq k \leq 2l-3$. By computing the action of $\dot{\tau}_i[k]$ on $\bd{s}_l$, one can show the following: \\
\begin{itemize}
\item[(i)] $s^{(k)}_{i_{k+1}}$, $s^{(k)}_{i_{k+1}+1} \in \set{s_1,\ldots,s_l,s_l+n,s_i+n,s_{i+1}-n}$. \\
\item[(ii)] If $s^{(k)}_{i_{k+1}}=s_a+\epsilon n$ and $s^{(k)}_{i_{k+1}+1}=s_{a'}+\epsilon' n$
with $a$, $a' \in \interv{1}{l}$, $\epsilon$, $\epsilon' \in \set{-1,0,1}$, then $a \neq a'$ and $\epsilon \epsilon'=0$. \\ 
\end{itemize}
Note moreover that by assumption on $\bd{s}_l$, we have for all $a$, $b \in \interv{1}{l}$ such that $a \neq b$,
$\big| s_b-s_a \big| \geq M \big| b-a \big| \geq M$. This together with (i), (ii) imply that
$$\big| s^{(k)}_{i_{k+1}} - s^{(k)}_{i_{k+1}+1} \big| \geq M-n.$$
As a consequence, $\Gamma(M-n)$ contains the arrow $\bd{s}_l^{(k)} \stackrel{i_k}{\longrightarrow} \bd{s}_l^{(k+1)}$ or 
$\bd{s}_l^{(k+1)} \stackrel{i_k}{\longrightarrow} \bd{s}_l^{(k)}$. It follows that $\bd{s}_l$ and $\bd{t}_l$ are in the same connected component
of $\Gamma(M-n)=\Gamma\big(n(N_0(w;\bd{r}_l)+1)\big)$. We can therefore apply Proposition \ref{prop_comparaison_bc_2} to conclude. 
\end{proof}

\begin{prop} \label{prop2_thm_comparaison_bc_cas_dominant}
Let $M \in \mN$. Then there exists $c \in \mZ$ (which only depends on $l$ and $n$, but not on $M$) such that
for all $(M+c)$-dominant multi-charges $\bd{s}_l$, $\bd{t}_l \in \mZ^l(s)$ with
$\mathcal{L}_{\bd{s}_l} = \mathcal{L}_{\bd{t}_l}$, we have $\bd{s}_l \Zconnected{M} \bd{t}_l$. 
\end{prop}

\begin{proof} We shall prove this proposition in Section \ref{pv_prop2_thm_comparaison_bc_cas_dominant}. \end{proof} 

\noindent \emph{Proof of Theorem \ref{thm_comparaison_bc_cas_dominant} from Propositions 
\ref{prop1_thm_comparaison_bc_cas_dominant} and \ref{prop2_thm_comparaison_bc_cas_dominant}.} \\
Let $M:=n(N_0(w;\bd{r}_l)+2)$. Let $c \in \mZ$ be the integer given by Proposition 
\ref{prop2_thm_comparaison_bc_cas_dominant} and put $N:=M+c$. Let $\bd{s}_l$, $\bd{t}_l \in \dot{W}_l.\bd{r}_l$ be two 
$N$-dominant multi-charges such that $\mathcal{L}_{\bd{s}_l}=\mathcal{L}_{\bd{t}_l}$. Put
 $\mathcal{L}:=\mathcal{L}_{\bd{s}_l}=\mathcal{L}_{\bd{t}_l}$. By Proposition 
 \ref{prop2_thm_comparaison_bc_cas_dominant}, there exist $\bd{s}_l^{(0)},\ldots,\bd{s}_l^{(r)} \in \mathcal{L} \cap \mathcal{C}_M$ such that 
$\bd{s}_l^{(0)}=\bd{s}_l$, $\bd{s}_l^{(r)}=\bd{t}_l$ and for all $1 \leq i \leq r$, we have 
$\bd{s}_l^{(i)} \in \set{\dot{\tau}_j^{\pm 1}.\bd{s}_l^{(i-1)} \mid 1 \leq j \leq l-1}$. Let $1 \leq i \leq r$. Since
$\mathcal{L}_{\bd{s}_l^{(i-1)}}=\mathcal{L}_{\bd{s}_l^{(i)}}=\mathcal{L}$, we have $\bd{s}_l^{(i-1)}$, 
$\bd{s}_l^{(i)} \in \dot{W}_l.\bd{r}_l$. By Proposition \ref{prop1_thm_comparaison_bc_cas_dominant}, the canonical 
bases of $\Fq[\bd{s}_l^{(i-1)}] \langle w+d(\bd{r}_l,\bd{s}_l^{(i-1)})\,\delta \rangle$ and $\Fq[\bd{s}_l^{(i)}] \langle w+d(\bd{r}_l,\bd{s}_l^{(i)})\,\delta \rangle$ are 
similar. Theorem \ref{thm_comparaison_bc_cas_dominant} follows. \cqfd

\subsection{Proof of Proposition \ref{prop2_thm_comparaison_bc_cas_dominant}} 
\label{pv_prop2_thm_comparaison_bc_cas_dominant}

The idea of the proof is the following. Let $\bd{s}_l$ and $\bd{t}_l$ be two multi-charges as in Proposition \ref{prop2_thm_comparaison_bc_cas_dominant} and put 
$\mathcal{L}:=\mathcal{L}_{\bd{s}_l} = \mathcal{L}_{\bd{t}_l}$. First, we introduce a suitable change of coordinates $\varphi$ that
maps $\mathcal{L}$ to $\mZ^{l-1}$ and such that $\bd{a}_l \Zconnected{M} \bd{b}_l$ if and only if $\varphi(\bd{a}_l)$ is $\mZ$-connected to
$\varphi(\bd{b}_l)$, that is there exists a piecewise affine path connecting $\varphi(\bd{a}_l)$ to $\varphi(\bd{b}_l)$ with edges parallel 
to the axes of coordinates of $\mZ^{l-1}$ (see Lemma \ref{lemme_phi_psi_cones}).
Roughly speaking, the aim is to replace the lattice $\mathcal{L}$ by $\mZ^{l-1}$ and the action of $\dot{\mathcal{Q}}$ by the obvious action
of $\mZ^{l-1}$ by translations. Doing this, we replace the set of $M$-dominance $\mathcal{C}_M$ by a certain cone (that is, an intersection of half-spaces),  
temporarily denoted by $C_M$. Note that two arbitrary points in $C_M \cap \mZ^{l-1}$ are not necessarily $\mZ$-connected in 
$C_M$. However, and this is the second step of the proof, we shall construct an integer $c$ such that $C_{M+c} \subset C_M$ and any
two points in $C_{M+c} \cap \mZ^{l-1}$ are $\mZ$-connected in $C_M$ (see Proposition \ref{prop_b_M_appartient_a_A_b_N}). 

\begin{notation} In addition to the notation from Section \ref{notation_paragraphe_comparaison_bc_cas_dominant} and
Proposition \ref{prop2_thm_comparaison_bc_cas_dominant}, we shall use for the proof the following notation.

\begin{itemize}
\item[*] For $\bd{x}=(x_1,\ldots,x_N) \in \mR^N$ ($N \in \mN^{*}$), put 
\begin{equation}
\lfloor \bd{x} \rfloor:=(\lfloor x_1 \rfloor,\ldots, \lfloor x_N \rfloor ) \in \mZ^N
\end{equation} 
and define $\lceil \bd{x} \rceil$ in a similar way. \\

\item[*] Let 
\begin{equation}
A:= \left( 
\begin{array}{cccccc}
2   & -1 & 0 & \ldots & 0	 \\
-1  & 2		 & -1 & \ddots  & \vdots \\
 0   & \ddots & \ddots & \ddots & 0 \\
\vdots& \ddots 		 & -1 	  & 2	   & -1 \\
0 &	\ldots & 0	  & -1	   & 2 \\
\end{array} \right) 
\end{equation}
denote the Cartan matrix of $\gsl{l}$. In particular, $A$ has $l-1$ rows and $l-1$ columns. \\

\item[*] For $1 \leq j \leq l-1$, let $\epsilon_j:=(\delta_{i,j})_{1 \leq i \leq l-1}$ be the $j$-th vector of the natural basis of $\mR^{l-1}$.
Put also $\bd{1}:=\epsilon_1+\cdots+\epsilon_{l-1}$. \\

\item[*] Define a partial ordering on the set of matrices (or vectors) of a given size with entries in $\mR$ by writing
$\bd{x}=(x_i)_{i \in I} \leq \bd{y}=(y_i)_{i \in I}$ if $x_i \leq y_i$ for all $i \in I$. We may speak of the maximum
of $\bd{x}$ and $\bd{y}$ (with respect to this ordering). Now, for $\bd{b} \in \mR^{l-1}$ define the cones
\begin{equation}
C_{\bd{b}} := \set{ \bd{x} \in \mR^{l-1} \mid A.\bd{x} \geq \bd{b}} \qquad \mbox{and} \qquad 
  C'_{\bd{b}} := \set{ \bd{x} \in \mR^{l-1} \mid A.\bd{x} \leq \bd{b}}.
\end{equation}  
The unique vector $\bd{\omega}:=\bd{\omega}(\bd{b}) \in \mR^{l-1}$ such that $A.\bd{\omega}=\bd{b}$ 
is called the \emph{vertex} of $C_{\bd{b}}$ (or $C'_{\bd{b}}$). \\

\item[*] For $M \in \mN$, let $\bd{b}(M)=\bd{b}(M;n,l,\bd{r_l})=(b_1^{(M)},\ldots,b_{l-1}^{(M)}) \in \mR^{l-1}$ denote the vector defined by
\begin{equation}
b_i^{(M)}:=(M+r_{i+1}-r_i)/n \qquad (1 \leq i \leq l-1).
\end{equation}
  
\item[*] Define a map 
$\varphi : (s_1,\ldots,s_l) \in \mR^l(s) \mapsto (x_1,\ldots,x_{l-1}) =  \varphi(s_1,\ldots,s_l) \in \mR^{l-1}$ by
\begin{equation}
x_i:= \frac{1}{n}{\sum_{j=1}^{i}(s_j-r_j)} \qquad (1 \leq i \leq l-1).
\end{equation}
Conversely, let $\psi : (x_1,\ldots,x_{l-1}) \in \mR^{l-1} \mapsto (s_1,\ldots,s_l)= \psi (x_1,\ldots,x_{l-1}) \in \mR^l(s)$
be the map defined by
\begin{equation}
s_i:=n(x_i-x_{i-1})+r_i \qquad (1 \leq i \leq l),
\end{equation}
where we put $x_0=x_l:=0$. \fini
\end{itemize} 
\end{notation}

\begin{lemma} \label{lemme_phi_psi_cones} \ \\[-5mm]
\begin{itemize}
\item[\rm (i)] The maps $\varphi : \mR^l(s) \rightarrow \mR^{l-1}$ and $\psi : \mR^{l-1} \rightarrow \mR^l(s)$ are bijections inverse to each 
other. \\
\item[\rm (ii)] We have $\psi(\mZ^{l-1})=\mathcal{L}_{\bd{r}_l}$, and for $1 \leq i \leq l-1$, 
$(x_1,\ldots,x_{l-1}) \in \mZ^{l-1}$, we have 
$$\psi(x_1,\ldots,x_i+1,\ldots,x_{l-1}) = \dot{\tau}_i. \psi(x_1,\ldots,x_{l-1}).$$ 
\item[\rm (iii)] For $M \in \mN$, we have $\varphi(\mathcal{C}_M)=C_{\bd{b}(M)}$.
\end{itemize}
\end{lemma}

\begin{proof} The proof of (i) and (ii) is straightforward. With obvious notation, we have the equivalence 
$$s_i-s_{i+1} \geq M \Longleftrightarrow -x_{i-1}+2x_i - x_{i+1} \geq (M+r_{i+1}-r_i)/n\,;$$
Statement (iii) follows.
\end{proof}

\begin{definition} Let $D \subset \mR^{l-1}$. We say that $\bd{x}$, $\bd{y} \in D$ are \emph{$\mZ$-connected in $D$} if there exist vectors 
$\bd{x}^{(0)},\ldots,\bd{x}^{(N)} \in D \cap \mZ^{l-1}$ such that $\bd{x}^{(0)}=\bd{x}$, $\bd{x}^{(N)}=\bd{y}$ and for all 
$0 \leq i \leq N-1$, we have $\bd{x}^{(i+1)}-\bd{x}^{(i)} \in \set { \pm \epsilon_j \mid 1 \leq j \leq l-1} ;$ in particular, we
have $\bd{x}$, $\bd{y} \in \mZ^{l-1}$. In this case, write $\bd{x} \zconnected{D} \bd{y}$. \fini
\end{definition}

In order to prove Proposition \ref{prop2_thm_comparaison_bc_cas_dominant}, we have to deal with $\mZ$-connected points in
cones $C_{\bd{b}}$ ($\bd{b} \in \mR^{l-1}$). Note that two points in $C_{\bd{b}} \cap \mZ^{l-1}$ are not necessary $\mZ$-connected in 
$C_{\bd{b}}$. For example, take $l=3$ (so $l-1=2$) and $\bd{b}=(0,0)$. Then $\bd{x}:=(1,1)$ and $\bd{y}:=(0,0)$ are two points in $C_{\bd{b}}$ 
which are not $\mZ$-connected in $C_{\bd{b}}$, because none of the points $(0,\pm 1)$ and $(\pm 1,0)$ lies in $C_{\bd{b}}$. However, given
$\bd{b} \in \mR^{l-1}$, one can construct $\bd{c} \leq \bd{b}$ such that any two points in $C_{\bd{b}}$ with integer coordinates are 
$\mZ$-connected in $C_{\bd{c}}$. This leads to the introduction of the following set. For $\bd{b} \in \mR^{l-1}$, put
\begin{equation}
\mathcal{A}(\bd{b}):=\set{\bd{c} \in \mR^{l-1} \mid \bd{c} \leq \bd{b} \mbox{ and } 
\forall \, \bd{x},\bd{y} \in C_{\bd{b}},\, \bd{x} \zconnected{C_{\bd{c}}} \bd{y}}.
\end{equation}
 
\begin{prop} \label{prop_A_b} Let $\bd{b} \in \mR^{l-1}$. Then we have the following.
\begin{itemize}
\item[\rm (i)] The set $\mathcal{A}(\bd{b})$ is nonempty.
\item[\rm (ii)] For all $\bd{c} \in \mathcal{A}(\bd{b})$, we have $\bd{c'} \leq \bd{c} \Rightarrow \bd{c'} \in \mathcal{A}(\bd{b})$.
\item[\rm (iii)] The map $\bd{b} \mapsto \mathcal{A}(\bd{b})$ is increasing.
\item[\rm (iv)] For all $\bd{b'} \in \mR^{l-1}$ such that $A^{-1}.(\bd{b}-\bd{b'}) \in \mZ^{l-1}$, we have 
$\mathcal{A}(\bd{b})=\mathcal{A}(\bd{b'})+(\bd{b}-\bd{b'})$.
\end{itemize}
\end{prop}

\begin{proof} We prove (i) and leave the other statements to the reader. Let $\bd{\omega}=(\omega_1,\ldots,\omega_{l-1})$ be the vertex of $C_{\bd{b}}$, 
$\bd{b'}:=\bd{b}+2.\bd{1}$ and $\bd{\omega}'$ be the vertex of $C_{\bd{b'}}$. Let
$\bd{\omega}''=(\omega''_1,\ldots,\omega''_{l-1})$ be equal to $\max(\bd{\omega}',\lceil \bd{\omega} \rceil)$ and let $\omega''_0=\omega''_l:=0$.
Now let $\bd{c}=(c_1,\ldots,c_{l-1}) \in \mR^{l-1}$ be the vector defined by
$$c_i:=-\omega''_{i-1}+2 \omega_i-\omega''_{i+1} \qquad (1 \leq i \leq l-1).$$  
We shall show that $\bd{c} \in \mathcal{A}(\bd{b})$. If is well-known that $A^{-1} \geq 0$, whence $\bd{x} \in C_{\bd{b}} \Rightarrow \bd{x} \geq \bd{\omega}$ and
$\bd{y} \in C'_{\bd{b'}} \Rightarrow \bd{y} \leq \bd{\omega}'$. Since $\bd{b} \leq \bd{b'}$, we get $\bd{\omega} \in C'_{\bd{b'}}$ and
therefore $\bd{\omega} \leq \bd{\omega}' \leq \bd{\omega}''$. Consider the set    
$$\textbf{P}:=\set{\bd{x} \in \mR^{l-1} \mid \bd{\omega} \leq \bd{x} \leq \bd{\omega}''}.$$
The argument above shows that $C_{\bd{b}} \cap C'_{\bd{b'}} \subset \textbf{P}$, and by construction we also have 
$\lceil \bd{\omega} \rceil \in \textbf{P}$.  By definition of $\bd{c}$, we have 
$\bd{\omega} \in \textbf{P} \subset C_{\bd{c}}$, whence $\bd{c} \leq \bd{b}$. Let us now show that any 
$\bd{x} \in C_{\bd{b}} \cap \mZ^{l-1}$ is $\mZ$-connected to
$\lceil \bd{\omega} \rceil$ in $C_{\bd{c}}$. Note that for 
$\bd{x}=(x_1,\ldots,x_{l-1}) \in C_{\bd{b}} \cap \mZ^{l-1}$, we have $\bd{x} \geq \bd{\omega}$ (because $\bd{x} \in C_{\bd{b}}$) and therefore 
$\bd{x} \geq \lceil \bd{\omega} \rceil$. We can thus argue by induction on 
$$N(\bd{x}):=\ds \sum_{i=1}^{l-1}(x_i-\lceil \omega_i \rceil) \in \mN.$$
If $N(\bd{x})=0$, we have $\bd{x}=\lceil \bd{\omega} \rceil \in \textbf{P} \subset C_{\bd{c}}$ and we are done. Assume now that 
$\bd{x} \in C_{\bd{b}} \cap \mZ^{l-1}$, $N(\bd{x})>0$, and consider two cases. Assume first that $\bd{x} \in C'_{\bd{b'}}$. Then we have 
$\bd{x} \in C_{\bd{b}} \cap C'_{\bd{b'}} \subset \textbf{P}$; moreover, we have 
$\lceil \bd{\omega} \rceil \in \textbf{P}$, so $\bd{x} \zconnected{\textbf{P}} \lceil \bd{\omega} \rceil$ because any two points
in $\textbf{P} \cap \mZ^{l-1}$ are $\mZ$-connected. Since $\textbf{P} \subset C_{\bd{c}}$, we can conclude in this case. Assume now that 
$\bd{x} \notin C'_{\bd{b'}}$. Let $1 \leq i \leq l-1$ be such that 
$-x_{i-1} + 2 x_i - x_{i+1} > b_i+2$ (where we put $x_0=x_l:=0$). Consider the vector
$$\bd{y}=(y_1,\ldots,y_{l-1}):=(x_1,\ldots,x_{i-1},x_i-1,x_{i+1},\ldots,x_{l-1}) \in \mZ^{l-1}$$
and put $y_0=y_l:=0$. Note that for $1 \leq j \leq l-1$, we have $$-y_{j-1} + 2 y_j - y_{j+1} \geq -x_{j-1} + 2 x_j - x_{j+1} - 2 \delta_{i,j}.$$
Using this fact together with the definition of $i$ and the assumption $\bd{x} \in C_{\bd{b}}$, we get $\bd{y} \in C_{\bd{b}}$. Moreover, since
$N(\bd{y})=N(\bd{x})-1$, we have by induction $\bd{y} \zconnected{C_{\bd{c}}} \lceil \bd{\omega} \rceil$. Note also that
$\bd{x} \zconnected{C_{\bd{c}}} \bd{y}$ because $\bd{c} \leq \bd{b}$, hence $\bd{x} \zconnected{C_{\bd{c}}} \lceil \bd{\omega} \rceil$. 
\end{proof}

Given $M \in \mN$, we shall construct $c \in \mZ$ such that $\bd{b}(M) \in \mathcal{A} \bigl( \bd{b}(M+c) \bigr)$.
Let $\bd{b} \in \mR^{l-1}$. By Proposition \ref{prop_A_b}, we can define
\begin{equation} \label{notation_suite_m_M}
m_M=m_M(n,l,\bd{r}_l):= \ds \max_{\bd{c}}{\min_{1 \leq i \leq l-1}{(c_i)}} \in \mZ,
\end{equation}
where $\bd{c}=(c_1,\ldots,c_{l-1})$ ranges over the set $\mathcal{A}\big(\bd{b}(M;n,l,\bd{r}_l)\big) \cap \mZ^{l-1}$.
\begin{lemma} \label{lemme_suite_m_M}
The sequence $(m_M)_{M \in \mN}$ is increasing. Moreover, we have $m_{nlM} = lM+m_0$ for all $M \in \mN$.
\end{lemma}

\begin{proof} The first statement follows from the fact that the maps $M \rightarrow \bd{b}(M)$ and $\bd{b} \mapsto \mathcal{A}(\bd{b})$ are increasing
(the latter by Proposition \ref{prop_A_b}). Let $M \in \mN$, $\bd{b}:=\bd{b}(nlM)$ and $\bd{b'}:=\bd{b}(0)$. Since $\det(A)=l$, the
Cramer formula shows that $A^{-1}.(\bd{b}-\bd{b'}) \in \mZ^{l-1}$. Applying Proposition \ref{prop_A_b}~(iv) to the pair $(\bd{b},\bd{b'})$
yields $\mathcal{A}\big( \bd{b}(nlM) \big)=\mathcal{A}\big( \bd{b}(0) \big)+(lM,\ldots,lM)$. The second statement follows.
\end{proof}

\begin{prop} \label{prop_b_M_appartient_a_A_b_N} 
Denote by $m_0^{(\min)}$ the real number $\min_{\bd{a}_l \in A_{l,n}(s)}m_0(n,l;\bd{a}_l)$ and define 
$c:=n+\lceil-m_0^{(\min)} n \rceil+nl \in \mZ$. Then for all $M \in \mN$, we have 
$\bd{b}(M) \in \mathcal{A} \bigl( \bd{b}(M+c) \bigr)$.
\end{prop}

\begin{proof} Let $\bd{c}=(c_1,\ldots,c_{l-1}) \in \mathcal{A}(\bd{b}(M+c)) \cap \mZ^{l-1}$ be such that 
$\min_{1 \leq i \leq l-1}(c_i)=m_{M+c}$. By proposition \ref{prop_A_b} (ii), it is enough to show that $\bd{b}(M) \leq \bd{c}$. We
have
$$M+c \geq nl \left(\ds \frac{M+n-m_0n}{nl}+1 \right) \geq nl \left \lceil \ds \frac{M+n-m_0n}{nl} \right \rceil 
= nl \left \lceil \ds \frac{\big(\frac{M}{n}+1 \big)-m_0}{l} \right \rceil.$$ 
The previous lemma implies 
$$m_{M+c} \geq l \left \lceil \ds \frac{\big(\frac{M}{n}+1 \big)-m_0}{l} \right \rceil + m_0 
    \geq l \left( \ds \frac{\big(\frac{M}{n}+1 \big)-m_0}{l} \right) + m_0 \geq \left \lceil \ds \frac{M}{n} \right \rceil. 
$$ 
Let $1 \leq i \leq l-1$. Since $\bd{r}_l \in A_{l,n}(s)$, we have $r_{i+1}-r_i \leq 0$, whence
$$\ds c_i \geq m_{M+c} \geq \left \lceil \frac{M}{n} \right \rceil \geq \frac{M}{n} + \frac{r_{i+1}-r_i}{n} =b_i^{(M)}.$$
The result follows.
\end{proof}

\begin{remark} \label{rq_constante_prop2_thm_comparaison_bc_cas_dominant}
The author proved in \cite[Section 4.4.3]{Y1} that $m_0^{(\min)} \geq -l^2$, which gives an explicit lower
bound for the integer $c$ from Propositions \ref{prop_b_M_appartient_a_A_b_N} and \ref{prop2_thm_comparaison_bc_cas_dominant}.\fini
\end{remark} 

\noindent \emph{Proof of Proposition \ref{prop2_thm_comparaison_bc_cas_dominant}. } \\
Let $c \in \mZ$ be the integer defined in Proposition \ref{prop_b_M_appartient_a_A_b_N}.
 Let $\bd{s}_l$, $\bd{t}_l \in \mZ^l(s)$ be two $(M+c)$-dominant multi-charges such that $\mathcal{L}_{\bd{s}_l} = \mathcal{L}_{\bd{t}_l}$. 
Let $\bd{r}_l \in A_{l,n}(s)$ be such that $\mathcal{L}_{\bd{r}_l} = \mathcal{L}_{\bd{s}_l} = \mathcal{L}_{\bd{t}_l}$. By Lemma \ref{lemme_phi_psi_cones}, we have 
$\varphi(\bd{s}_l)$, $\varphi(\bd{t}_l) \in C_{\bd{b}(M+c)} \cap \mZ^{l-1}$. Proposition \ref{prop_b_M_appartient_a_A_b_N} now implies that
$\bd{b}(M)$ is in $\mathcal{A} \bigl( \bd{b}(M+c) \bigr)$, whence $\varphi(\bd{s}_l) \zconnected{C_{\bd{b}(M)}} \varphi(\bd{t}_l)$. 
Applying again Lemma \ref{lemme_phi_psi_cones} yields $\bd{s}_l \Zconnected{M} \bd{t}_l$. \cqfd

\vspace{5mm}
\small \rm Xavier YVONNE, Laboratoire de Math\'ematiques Nicolas Oresme, Universit\'e de Caen, BP 5186, 14032 Caen Cedex, France. \\

\indent \it E-mail address: \tt xyvonne@math.unicaen.fr

\clearpage
\begin{thebibliography}{ABCD}

\bibitem[A]{A} \scshape S. Ariki, \rm \emph{On the decomposition numbers of the Hecke algebra of $G(m,1,n)$}, J. Math. Kyoto Univ. 
\textbf{36}, no. 4 (1996), 789-808. 

\bibitem[DJM]{DJM} \scshape R. Dipper, G. James, A. Mathas, \rm \emph{Cyclotomic $q$-Schur algebras}, Math. Z. \textbf{229} (1998), 385-416.

\bibitem[FLOTW]{FLOTW} \scshape O. Foda, B. Leclerc, M. Okado, J.-Y. Thibon, T. Welsh, \rm \emph{Branching functions of $A_{n-1}^{(1)}$ 
and Jantzen-Seitz problem for Ariki-Koike algebras}, Adv. Math. \textbf{141} No.2 (1999), 322-365. 

\bibitem[H]{H} \scshape T. Hayashi, \rm \emph{$q$-analogues of Clifford and Weyl algebras - spinor and oscillator representations of quantum 
enveloping algebras}, Commun. Math. Phys. \textbf{127} No.1 (1990), 129-144. 

\bibitem[J]{J} \scshape N. Jacon, \rm \emph{On the parametrization of the simple modules for Ariki-Koike algebras at roots of unity}, 
J. Math. Kyoto Univ. \textbf{44} No.4 (2004), 729-767.
 
\bibitem[JMMO]{JMMO} \scshape M. Jimbo, K. Misra, T. Miwa, M. Okado, \rm \emph{Combinatorics of representations of $\Uq$ at $q=0$}, 
Commun. Math. Phys. \textbf{136} No.3 (1991), 543-566.

\bibitem[Kac]{Kac} \scshape V. G. Kac, \rm \emph{Infinite dimensional Lie algebras}, 3rd Ed. Cambridge University Press, 1990.

\bibitem[Kas1]{Kas1} \scshape M. Kashiwara, \rm \emph{Global crystal bases of quantum groups}, Duke Math. J. \textbf{69} (1993), 455-485.

\bibitem[Kas2]{Kas2} \scshape M. Kashiwara, \rm \emph{The crystal base and Littelmann's refined Demazure character formula}, Duke Math. J. 
\textbf{71} (1993), 839-958.

\bibitem[KMS]{KMS} \scshape M. Kashiwara, T. Miwa, E. Stern, \rm \emph{Decomposition of $q$-deformed Fock spaces}, 
Selecta Math. \textbf{1} (1995), 787-805.

\bibitem[KT]{KT} \scshape M. Kashiwara, T. Tanisaki, \rm \emph{Parabolic Kazhdan-Lusztig polynomials and Schubert varieties}, J. Algebra 
\textbf{249} (2002), 306-325.

\bibitem[LLT]{LLT} \scshape A. Lascoux, B. Leclerc, J.-Y. Thibon, \rm \emph{Hecke algebras at roots of unity and crystal bases
of quantum affine algebras}, Commun. Math. Phys. \textbf{181} No.1 (1996), 205-263.

\bibitem[LM]{LM} \scshape B. Leclerc, H. Miyachi, \rm \emph{Some closed formulas for canonical bases of Fock spaces},
Represent. Theory \textbf{6} (2002), 290-312.

\bibitem[LT1]{LT1} \scshape B. Leclerc, J.-Y. Thibon, \rm \emph{Canonical bases of $q$-deformed Fock spaces}, Internat. Math. Res.
Notices \textbf{9} (1996), 447-456.

\bibitem[LT2]{LT2} \scshape B. Leclerc, J.-Y. Thibon, \rm \emph{Littlewood-Richardson coefficients and Kazhdan-Lusztig polynomials}, 
Combinatorial Methods in Representation Theory, Advanced Studies in Pure Mathematics \textbf{28} (2000), 155-220.

\bibitem[Mac]{Mac} \scshape I. G. Macdonald, \rm \emph{Symmetric functions and Hall polynomials}, 2nd Ed. Oxford Science Publications, 
Oxford University Press, 1995.

\bibitem[MM]{MM} \scshape K.C. Misra, T. Miwa, \rm \emph{Crystal base for the basic representation of $\Uq$},
Commun. Math. Phys. \textbf{134} No.1 (1990), 79-88.

\bibitem[S]{S} \scshape J. Scopes, \rm \emph{Cartan matrices and Morita equivalence for blocks of the symmetric groups},
J. Algebra \textbf{142} No.2 (1991), 441-455.

\bibitem[U]{U} \scshape D. Uglov, \rm \emph{Canonical bases of higher-level $q$-deformed Fock spaces and Kazhdan-Lusztig polynomials}, 
in Physical Combinatorics ed. M. Kashiwara, T. Miwa, Progress in Math. \textbf{191}, Birkhäuser (2000), arXiv math.QA/9905196 (1999).  

\bibitem[VV]{VV} \scshape M. Varagnolo, E. Vasserot, \rm \emph{On the decomposition matrices of the quantized Schur algebra}, Duke Math. J. 
\textbf{100} (1999), 267-297.

\bibitem[Y1]{Y1} \scshape X. Yvonne, \rm \emph{Bases canoniques d'espaces de Fock de niveau sup\'erieur}, Thèse de l'Universit\'e de Caen (2005).

\bibitem[Y2]{Y2} \scshape X. Yvonne, \rm \emph{A conjecture for $q$-decomposition matrices of cyclotomic $v$-Schur algebras}, 
J. Algebra (2006), arXiv math.RT/0505379. 
\end{thebibliography}
\end{document}